

\input amssym  

\newdimen\normalparindent

\iffalse  

\hsize=15truecm
\hoffset=.46truecm
\vsize=23.7truecm
\voffset=.46truecm

\normalparindent=24pt

\font\elevenrm=cmr10 at 11pt 
\font\eightrm=cmr8 
\font\sixrm=cmr6 

\font\eleveni=cmmi10 at 11pt 
\font\eighti=cmmi8
\font\sixi=cmmi6

\font\elevensy=cmsy10 at 11pt 
\font\eightsy=cmsy8
\font\sixsy=cmsy6

\font\elevenex=cmex10 at 11pt 

\font\elevenbf=cmbx10 at 11pt 
\font\eightbf=cmbx8
\font\sixbf=cmbx6

\font\eleventt=cmtt10 at 11pt 
\font\elevensl=cmsl10 at 11pt 
\font\elevenit=cmti10 at 11pt 

\textfont0=\elevenrm \scriptfont0=\eightrm \scriptscriptfont0=\sixrm
\def\rm{\fam0\elevenrm}
\textfont1=\eleveni \scriptfont1=\eighti \scriptscriptfont1=\sixi
 
\textfont2=\elevensy \scriptfont2=\eightsy \scriptscriptfont2=\sixsy

\textfont3=\elevenex \scriptfont3=\elevenex \scriptscriptfont3=\elevenex
\textfont\itfam=\elevenit
\def\it{\fam\itfam\elevenit}
\textfont\slfam=\elevensl
\def\sl{\fam\slfam\elevensl}
\textfont\bffam=\elevenbf \scriptfont\bffam=\eightbf
\scriptscriptfont\bffam=\sixbf
\def\bf{\fam\bffam\elevenbf}
\textfont\ttfam=\eleventt
\def\tt{\fam\ttfam\eleventt}

\skewchar\eleveni='177 \skewchar\eighti='177 \skewchar\sixi='177
\skewchar\elevensy='60 \skewchar\eightsy='60 \skewchar\sixsy='60

\font\sc=cmcsc10 at 11pt 

\font\cmcyr=cmcyr10 at 11pt
\font\cmcti=cmcti10 at 11pt
\font\cmccsc=cmccsc10 at 11pt

\smallskipamount=3.5pt plus 1pt minus 1pt
\medskipamount=7pt plus 2pt minus 2pt
\bigskipamount=14pt plus 2pt minus 2pt
\normalbaselineskip=14pt
\normallineskip=1pt
\normallineskiplimit=0pt
\jot=3.5pt

\normalbaselines
\rm

\else  

\hsize=15truecm
\hoffset=.46truecm
\vsize=23.7truecm
\voffset=.46truecm

\normalparindent=20pt

\font\sc=cmcsc10

\font\cmcyr=cmcyr10
\font\cmcti=cmcti10
\font\cmccsc=cmccsc10

\fi

\def\cyrchardefs{%
\chardef\yo=60
\chardef\Yo=62
\chardef\yu=64
\chardef\a=65
\chardef\b=66
\chardef\ts=67
\chardef\d=68
\chardef\ye=69
\chardef\f=70
\chardef\g=71
\chardef\kh=72
\chardef\i=73
\chardef\j=74
\chardef\k=75
\chardef\l=76
\chardef\m=77
\chardef\n=78
\chardef\o=79
\chardef\p=80
\chardef\ya=81
\chardef\r=82
\chardef\s=83
\chardef\t=84
\chardef\u=85
\chardef\zh=86
\chardef\v=87
\chardef\soft=88
\chardef\y=89
\chardef\z=90
\chardef\sh=91
\chardef\e=92
\chardef\shch=93
\chardef\ch=94
\chardef\hard=95
\chardef\Yu=96
\chardef\A=97
\chardef\B=98
\chardef\Ts=99
\chardef\D=100
\chardef\Ye=101
\chardef\F=102
\chardef\G=103
\chardef\H=104
\chardef\I=105
\chardef\J=106
\chardef\K=107
\chardef\L=108
\chardef\M=109
\chardef\N=110
\chardef\O=111
\chardef\P=112
\chardef\Ya=113
\chardef\R=114
\chardef\S=115
\chardef\T=116
\chardef\V=119
\chardef\Soft=120
\chardef\Y=121
\chardef\Z=122
\chardef\Sh=123
\chardef\E=124
\chardef\Shch=125
\chardef\Ch=126
\chardef\Hard=127
}
\def\cyr{\cmcyr\cyrchardefs}
\def\cycsc{\cmccsc\cyrchardefs}
\def\cyti{\cmcti\cyrchardefs}

\def\S{\mathhexbox278\thinspace}

\def\square{\hbox to.77778em{%
\hfil\vrule\vbox to.675em{\hrule width.6em\vfil\hrule}\vrule\hfil}}

\long\def\acknowledgements#1\par{\medbreak\noindent{\it
    Acknowledgements\/}.\enspace #1\par\medbreak}
\def\definition#1\par{\medbreak\noindent{\bf Definition.}\enspace
  #1\par\medbreak}
\def\example#1\par{\medbreak\noindent{\bf Example.}\enspace
  #1\par\medbreak}
\long\def\remark#1\par{\medbreak\noindent{\it Remark\/}.\enspace
#1\par\medbreak}
\long\def\remarks#1\par{\medbreak\noindent{\it Remarks\/}.\enspace
#1\par\medbreak}
\def\exercise#1\par{\medbreak\noindent{\bf Exercise.}\enspace
#1\par\medbreak}
\def\notation#1\par{\medbreak\noindent{\bf Notation.}\enspace
#1\par\medbreak}
\def\proof{\noindent{\it Proof\/}.\enspace}
\def\endproof{\nobreak\hfill\quad\square\par\medbreak}

\def\lineover#1{{\offinterlineskip\mathchoice
{\setbox0=\hbox{$\displaystyle#1$}%
\vbox{\kern .33pt\hbox to\wd0{\kern 1pt\leaders\hrule height .33pt%
\hfill\kern 1pt}\kern 1pt\box0}}
{\setbox0=\hbox{$\textstyle#1$}%
\vbox{\kern .33pt\hbox to\wd0{\kern 1pt\leaders\hrule height .33pt%
\hfill\kern 1pt}\kern 1pt\box0}}
{\setbox0=\hbox{$\scriptstyle#1$}%
\vbox{\kern .25pt\hbox to\wd0{\kern .8pt\leaders\hrule height .25pt%
\hfill\kern .8pt}\kern .8pt\box0}}
{\setbox0=\hbox{$\scriptscriptstyle#1$}%
\vbox{\kern .2pt\hbox to\wd0{\kern .6pt\leaders\hrule height .2pt%
\hfill\kern .6pt}\kern .6pt\box0}}}}

\def\isom{\buildrel\sim\over\longrightarrow}
\def\morphism#1{\buildrel#1\over\longrightarrow}
\def\isomorphism#1{\mathrel{\mathop{\longrightarrow}%
\limits^{#1}_{\raise0.5ex\hbox{$\scriptstyle\sim$}}}}

\def\injlim{\mathop{\vtop{\offinterlineskip\halign{##\cr
 \hfil\rm lim\hfil\cr\noalign{\kern.1ex}\rightarrowfill\cr
 \noalign{\kern-.4ex}\cr}}}}
\def\projlim{\mathop{\vtop{\offinterlineskip\halign{##\cr
 \hfil\rm lim\hfil\cr\noalign{\kern.1ex}\leftarrowfill\cr
 \noalign{\kern-.4ex}\cr}}}}

\def\smallmatrix#1#2#3#4{\bigl({#1\atop#3}\,{#2\atop#4}\bigr)}
\def\medmatrix#1#2#3#4{\biggl({#1\atop#3}\,\,\,{#2\atop#4}\biggr)}

\def\blank{\mkern12mu}

\def\textfrac#1/#2{{\textstyle{#1\over#2}}}

\def\relativediag#1#2#3#4#5#6{\vcenter{\baselineskip=3ex \halign{
\hfil$##$&$##$&$##$\hfil\cr
#1\quad& \hfilneg\buildrel#2\over\longrightarrow\hfilneg& \quad#3\cr
\lower1ex\llap{$\scriptstyle#4\hskip-1ex$}\searrow& \quad&
\swarrow\lower1ex\rlap{$\hskip-1ex\scriptstyle#5$} \cr
& \hfilneg#6\hfilneg&\cr}}}

\def\trianglediag#1#2#3#4#5#6{\vcenter{\baselineskip=3ex \halign{
\hfil$##$&$##$&$##$\hfil\cr
#1\quad& \hfilneg\buildrel#2\over\longrightarrow\hfilneg& \quad#3\cr
\lower1ex\llap{$\scriptstyle#6\hskip-1ex$}\nwarrow& \quad&
\swarrow\lower1ex\rlap{$\hskip-1ex\scriptstyle#4$} \cr
& \hfilneg#5\hfilneg&\cr}}}

\def\correspondence#1#2#3#4#5{\vcenter{\baselineskip=3ex \halign{
\hfil$##$&$##$&$##$\hfil\cr
&\hfilneg#1\hfilneg\cr
\raise1ex\llap{$\scriptstyle#2$}\swarrow&&\searrow
\raise1ex\rlap{$\scriptstyle#3$}\cr
#4&&#5\cr}}}

\newif\iffirstpar
\everypar{\iffirstpar\parindent=\normalparindent\firstparfalse\fi}

\def\sectionheading#1{\subcount=0 \subsectioncount=0 \eqcount=0
  \bigskip\vskip\parskip
  \leftline{\bf #1}\nobreak\smallskip\firstpartrue\parindent=0pt}

\newif\ifappendix \appendixfalse

\def\currentsection{\ifappendix A\else\number\sectioncount\fi}

\def\section#1\par{\advance\sectioncount by1%
  \edef\currentlabel{\currentsection}%
  \sectionheading{\currentsection.\enspace#1}}

\def\unnumberedsection#1\par{\sectionheading{#1}}

\def\subsection#1\par{\medbreak\penalty-200\advance\subsectioncount by1%
  \edef\currentlabel{\currentsection.\number\subsectioncount}%
  \leftline{\it\currentsection.\number\subsectioncount.\enspace#1}%
  \smallskip\parindent=0pt\firstpartrue}

\newwrite\auxfile

\newcount\sectioncount \sectioncount=0
\newcount\subsectioncount 
\newcount\subcount 
\newcount\eqcount 

\def\subno{\global\advance\subcount by1\relax
  \currentsection.\number\subcount
  \xdef\currentlabel{\currentsection.\number\subcount}}
\def\proclaim #1. #2\par{\medbreak
  \noindent{\bf#1~\subno.\enspace}{\sl#2\par}%
  \ifdim\lastskip<\medskipamount \removelastskip\penalty55\medskip\fi}
\def\proclaimx #1 (#2). #3\par{\medbreak
  \noindent{\bf#1~\subno\ \rm (#2).\enspace}{\sl#3\par}%
  \ifdim\lastskip<\medskipamount \removelastskip\penalty55\medskip\fi}

\newdimen\algindent
\def\plusindent{\advance\algindent by \parindent}
\def\minusindent{\advance\algindent by-\parindent}


\newcount\algstepcount

\long\def\algorithm (#1). #2\endalgorithm{\medbreak
  \algindent=0pt%
  \algstepcount=0%
  \noindent{\bf Algorithm~\subno} (#1). {\sl#2}\par\medbreak}

\def\step{\advance\algstepcount by1
\edef\currentlabel{\number\algstepcount}
\smallskip\hangindent\parindent
\advance\hangindent by\algindent\indent
\llap{{\bf \the\algstepcount.}\enspace}\kern\algindent
\ignorespaces}


\def\labeldef#1#2{\expandafter\gdef\csname L@#1\endcsname{#2}}
\def\label#1{%
  \expandafter\xdef\csname L@#1\endcsname{\currentlabel}%
  \write\auxfile{\string\labeldef{#1}{\csname L@#1\endcsname}}%
  \ignorespaces}
\def\ref#1{{\rm\expandafter\ifx\csname L@#1\endcsname\relax
  \message{Undefined label `#1'}??\else
  \csname L@#1\endcsname\fi}}

\def\eqdef#1#2{\expandafter\gdef\csname E@#1\endcsname{#2}}
\def\eqnumber#1{\global\advance\eqcount by1\relax
  \eqno{\rm(\currentsection.\number\eqcount)}%
  \expandafter\xdef\csname E@#1\endcsname{%
    \currentsection.\number\eqcount}%
  \write\auxfile{\string\eqdef{#1}{\csname E@#1\endcsname}}}
\def\eqref#1{{\rm(\expandafter\ifx\csname E@#1\endcsname\relax
  \message{Undefined equation `#1'}??\else
  \csname E@#1\endcsname\fi)}}

\newcount\refcount \refcount=0
\def\citedef#1#2{\expandafter\gdef\csname C@#1\endcsname{#2}}
\def\cite#1{\expandafter\ifx\csname C@#1\endcsname\relax
  \message{Undefined reference `#1'}\citedef{#1}{??}\fi
  \expandafter\gdef\csname R@#1\endcsname{\relax}%
  [\csname C@#1\endcsname]}
\def\citex#1#2{\expandafter\ifx\csname C@#1\endcsname\relax
  \message{Undefined reference `#1'}\citedef{#1}{??}\fi
  \expandafter\gdef\csname R@#1\endcsname{\relax}%
  [\csname C@#1\endcsname, #2]}
\def\reference#1{\advance\refcount by 1%
  \expandafter\ifx\csname R@#1\endcsname\relax
  \message{Warning: reference `#1' not used}\fi
  \expandafter\edef\csname C@#1\endcsname{\the\refcount}%
  \write\auxfile{\string\citedef{#1}{\csname C@#1\endcsname}}%
  \item{[\csname C@#1\endcsname]}}

\newif\ifauxexists
\immediate\openin0=\jobname.aux
\ifeof 0
  \auxexistsfalse
\else
  \auxexiststrue
\fi
\immediate\closein0
\ifauxexists
  \input \jobname.aux
\else
  \message{No file `\jobname.aux'}
\fi
\openout\auxfile=\jobname.aux

\def\c{{\frak c}}
\def\C{{\bf C}}
\def\fd{{\frak d}}
\def\gr{\mathop{\rm gr}\nolimits}
\def\hyp{{\bf H}}
\def\Ltwo{{\rm L}^2}
\def\muhyp{\mu_\hyp^{}}
\def\P{{\bf P}}
\def\Q{{\bf Q}}
\def\R{{\bf R}}
\def\si{\mathop{\rm si}\nolimits}
\def\SL{{\rm SL}}
\def\vol{\mathop{\rm vol}\nolimits}
\def\Z{{\bf Z}}

\def\implies{\;\Longrightarrow\;}

\font\titlefont=cmssbx10 scaled \magstep2

\centerline{\titlefont An approximate spectral representation
and explicit bounds}
\smallskip
\centerline{\titlefont for Green functions of Fuchsian
groups}

\medskip
\centerline{Peter Bruin}
\smallskip
\centerline{19 July 2012}

\bigskip {\narrower\noindent {\it Abstract\/.}\enspace We study the
  Green function~$\gr_\Gamma$ for the Laplace operator on the quotient
  of the hyperbolic plane by a cofinite Fuchsian group~$\Gamma$.  We
  use a limiting procedure, starting from the resolvent kernel, and
  lattice point estimates for the action of~$\Gamma$ on the hyperbolic
  plane to prove an ``approximate spectral representation''
  for~$\gr_\Gamma$.  Combining this with bounds on Maa\ss\ forms and
  Eisenstein series for~$\Gamma$, we prove explicit bounds
  on~$\gr_\Gamma$.\par}

\vfootnote{}{{\it Mathematics Subject Classification (2010):\/}
11F72,
30F35, 35J08

This paper evolved from part of the author's thesis \cite{thesis}, the
research for which was supported by the Netherlands Organisation for
Scientific Research (NWO).  Further research was supported by Swiss
National Science Foundation grant 124737 and by the
Max-Planck-Institut f\"ur Mathematik, Bonn.}

\section Introduction and statement of results

\label{sec:introduction}

The hyperbolic plane~$\hyp$ is the unique two-dimensional, complete,
connected and simply connected Riemannian manifold with constant
Gaussian curvature~$-1$.  We identify $\hyp$ with the complex upper
half-plane; this gives $\hyp$ a complex structure.  In terms of the
standard coordinate~$z=x+iy$, the Riemannian metric is
$$
{dz\,d\bar z\over(\Im z)^2}={dx^2+dy^2\over y^2},
$$
and the associated volume form is
$$
\muhyp={i\,dz\wedge d\bar z\over 2(\Im z)^2}
={dx\wedge dy\over y^2}.
$$
Instead of using the geodesic distance $r(z,w)$ on~$\hyp$ directly, we
use the more convenient function
$$
\eqalign{
u(z,w)&=\cosh r(z,w)\cr
&=1+{|z-w|^2\over2(\Im z)(\Im w)}.}
$$

Let $\Delta$ denote the Laplace--Beltrami operator on~$\hyp$, given by
$$
\Delta=y^2(\partial_x^2+\partial_y^2).
$$
The Green function for~$\Delta$ is the unique smooth real-valued
function $\gr_\hyp$ outside the diagonal on~$\hyp\times\hyp$
satisfying
$$
\displaylines{
\gr_\hyp(z,w)={1\over2\pi}\log|z-w|+O(1)\quad\hbox{as }z\to w,\cr
\Delta\gr_\hyp(\blank,w)=\delta_w\quad\hbox{for all }w\in\hyp,\cr
\gr_\hyp(z,w)=O(u(z,w)^{-1})\quad\hbox{as }u(z,w)\to\infty,}
$$
where $\Delta$ is taken with respect to the first variable.  It is
given by
$$
\gr_\hyp(z,w)=-L(u(z,w)),
$$
where
$$
L(u)={1\over4\pi}\log{u+1\over u-1}.
\eqnumber{definition-L}
$$

The group~$\SL_2(\R)$ acts on~$\hyp$ by isometries.  Under the
identification of~$\hyp$ with the complex upper half-plane, this
action on~$\hyp$ is the restriction of the action on~$\P^1(\C)$ by
M\"obius transformations.  Elements of~$\SL_2(\R)\setminus\{\pm1\}$
are classified according to their fixed points in~$\P^1(\C)$ as
elliptic (two conjugate fixed points in $\P^1(\C)\setminus\P^1(\R)$),
parabolic (a unique fixed point in $\P^1(\R)$), and hyperbolic (two
distinct fixed points in $\P^1(\R)$).  This terminology also applies
to conjugacy classes.

A {\it Fuchsian group\/} is a discrete subgroup of~$\SL_2(\R)$.  A
{\it cofinite Fuchsian group\/} is a Fuchsian group~$\Gamma$ such that
the volume of~$\Gamma\backslash\hyp$ with respect to the measure
induced by~$\muhyp$ is finite.  We will exclusively consider cofinite
Fuchsian groups, and for such a group~$\Gamma$ we write
$$
\vol_\Gamma=\int_{\Gamma\backslash\hyp}\muhyp.
$$

\remark We define integration on $\Gamma\backslash\hyp$ in a
stack-like way, so the above integral is $1/\#(\Gamma\cap\{\pm1\})$
times the integral over $\Gamma\backslash\hyp$ viewed as a Riemann
surface.  This implies that if $f$ is a $\Gamma$-invariant function
on~$\hyp$ and $\Gamma'$ is a subgroup of finite index in~$\Gamma$,
then
$$
\int_{\Gamma'\backslash\hyp}f\muhyp=(\Gamma:\Gamma')
\cdot\int_{\Gamma\backslash\hyp}f\muhyp.
$$
Furthermore, this definition justifies the method of ``unfolding'': if
$f$ is a smooth function with compact support on~$\hyp$ and $F$ is the
function on~$\Gamma\backslash\hyp$ defined by
$$
F(z)=\sum_{\gamma\in\Gamma}f(\gamma z),
$$
then
$$
\int_{\Gamma\backslash\hyp}F\muhyp=
\int_\hyp f\muhyp.
$$

Let $\Gamma$ be a cofinite Fuchsian group.  The restriction of the
Laplace operator~$\Delta$ to the space of smooth and bounded
$\Gamma$-invariant functions on~$\hyp$ can be extended to an
(unbounded, densely defined) self-adjoint operator on the Hilbert
space~$\Ltwo(\Gamma\backslash\hyp)$, which we denote
by~$\Delta_\Gamma$.

The operator~$\Delta_\Gamma$ is invertible on the orthogonal
complement of the constant functions in the following sense: there
exists a unique bounded self-adjoint operator~$G_\Gamma$
on~$\Ltwo(\Gamma\backslash\hyp,\muhyp)$ such that for all smooth and
bounded functions $f$ on~$\Gamma\backslash\hyp$ the function~$G_\Gamma
f$ satisfies
$$
\Delta_\Gamma G_\Gamma f=f-{1\over \vol_\Gamma}\int_{\Gamma\backslash\hyp}f\muhyp
\quad\hbox{and}\quad\int_{\Gamma\backslash\hyp}G_\Gamma f\muhyp=0.
$$
There exists a unique function $\gr_\Gamma$ on $\hyp\times\hyp$ that
is $\Gamma$-invariant in both variables separately, satisfies
$\gr_\Gamma(z,w)=\gr_\Gamma(w,z)$, is smooth except for logarithmic
singularities at points of the form $(z,\gamma z)$, and has the
property that if $f$ is a smooth and bounded $\Gamma$-invariant
function on~$\hyp$, then the function~$G_\Gamma f$ is given by
$$
G_\Gamma f(z)=\int_{w\in\Gamma\backslash\hyp}\gr_\Gamma(z,w)f(w)\muhyp(w).
$$
The function~$\gr_\Gamma$ is called the {\it Green function of the
Fuchsian group~$\Gamma$\/}.




In this paper, we study $\gr_\Gamma$ quantitatively, with the goal of
obtaining explicit upper and lower bounds.  One result that can be
stated without introducing too much notation is the following.

\proclaimx Theorem (corollary of Theorem~\ref{theorem:explicit}).  Let
$\Gamma_0$ be a cofinite Fuchsian group, let $Y_0$ be a compact subset
of\/~$\Gamma_0\backslash\hyp$, and let $\delta>1$ and~$\eta>0$ be real
numbers.  There exist real numbers $A$ and~$B$ such that the following
holds.  Let $\Gamma$ be a subgroup of finite index in~$\Gamma_0$ such
that all non-zero eigenvalues of~$-\Delta_\Gamma$ are at least $\eta$.
Then for all $z,w\in\hyp$ whose images in~$\Gamma_0\backslash\hyp$ lie
in~$Y_0$, we have
$$
A\le\gr_\Gamma(z,w)
+\sum_{\textstyle{\gamma\in\Gamma\atop u(z,\gamma w)\le\delta}}
\bigl(L(u(z,\gamma w))-L(\delta)\bigr)
\le B.
$$

\label{theorem:introduction}

Let us give a very brief overview of the article.  In
Section~\ref{sec:tools}, we collect known results about Fuchsian
groups.  Most importantly, we make use of the hyperbolic lattice point
problem and the techniques used to attack this problem.  In
Section~\ref{section:approximate-spectral-representation}, we use
these results, together with a construction of the Green function
involving the resolvent kernel, to ``sandwich'' the Green
function~$\gr_\Gamma$ (with the logarithmic singularity removed)
between two functions that, unlike $\gr_\Gamma$ itself, admit spectral
representations.  In Section~\ref{sec:explicit-bounds}, we bound these
functions in a way that lends itself to explicit evaluation.  As an
example, we find explicit constants $A$ and~$B$ as in
Theorem~\ref{theorem:introduction} in the case where
$\Gamma_0=\SL_2(\Z)$, $\Gamma\subseteq\Gamma_0$ is a congruence
subgroup, $\delta=2$ and $Y_0$ is the compact subset
of~$\Gamma_0\backslash\hyp$ corresponding to the points $z$ in the
standard fundamental domain of~$\SL_2(\Z)$ such that $\Im z\le 2$.  In
Section~\ref{section:extension-cusps}, we use the bounds given by
Theorem~\ref{theorem:introduction} to deduce bounds
on~$\gr_\Gamma(z,w)$ in the case where $Y$ is obtained by cutting out
discs around the cusps of~$\Gamma$ and where $z$, $w$ or both are in
such a disc.  Finally, a number of bounds on Legendre functions that
we will need have been collected in an appendix.

In a forthcoming paper, we will use the results in this article to
obtain explicit bounds on the {\it canonical Green function\/} of a
modular curve~$X$.  This function is defined similarly to~$\gr_\hyp$,
using the {\it canonical $(1,1)$-form\/} on~$X$ instead of~$\muhyp$.
It plays a fundamental role in Arakelov theory; see
Arakelov~\cite{Arakelov} and Faltings~\cite{Faltings: Calculus on
  arithmetic surfaces}.  Explicit bounds on canonical Green functions
of modular curves are relevant to the work of Edixhoven, Couveignes
et~al.~\cite{book} and the author~\cite{thesis} on computing
two-dimensional representations of the absolute Galois group of~$\Q$
that are associated to Hecke eigenforms over finite fields.

This article may be compared with earlier work of Jorgenson and Kramer
on bounding canonical and hyperbolic Green functions of compact
Riemann surfaces \citex{Jorgenson-Kramer: Green's
  functions}{especially Theorem~4.5}.  Jorgenson and Kramer consider
compact Riemann surfaces $X$ of genus at least~2, which can be
obtained as $X=\Gamma\backslash\hyp$ for a cofinite Fuchsian
group~$\Gamma$ without elliptic and parabolic elements.  They obtain
bounds on the hyperbolic Green function by comparing it to the heat
kernel on~$X$.  Our method, too, starts with comparing the Green
function with a kernel that can be obtained as a sum over elements
of~$\Gamma$, but the subsequent arguments are rather different.  Let
us note some of the differences.  First, we allow arbitrary cofinite
Fuchsian groups, which is the natural setting for modular curves.
Second, the procedure that we apply
in~\S\ref{construction-green-function} to construct the Green function
as a limit of a family of kernels $K_a$ for $a\to 1$ leads to bounds
that are independent of the specific family.  We take $K_a$ to be the
resolvent kernel with parameter $a\to 1$, but the heat kernel with
parameter $t\to\infty$ could have been used with the same result; see
\citex{thesis}{\S II.5.2}.  Finally, our bounds are much easier to
make explicit than those in \cite{Jorgenson-Kramer: Green's
  functions}; this is illustrated in~\S\ref{subsection:example}.

\acknowledgements Part of this paper was written during a stay at the
Max-Planck-Institut f\"ur Mathematik in Bonn; I am grateful for its
hospitality.  I thank Ariyan Javanpeykar for comments on an earlier
version.  The computations outlined in \S\ref{subsection:example} were
carried out using PARI/GP~\cite{pari}.

\section Tools

\label{sec:tools}

\subsection Cusps

\label{cusps}

Let $\Gamma$ be a cofinite Fuchsian group.  The cusps of~$\Gamma$
correspond to the conjugacy classes of non-trivial maximal parabolic
subgroups in~$\Gamma$.  Every such subgroup has a unique fixed point
in~$\P^1(\R)$.  For every cusp~$\c$ we choose a representative of the
corresponding conjugacy class and denote it by~$\Gamma_\c$.


Let $\c$ be a cusp of~$\Gamma$.  We fix an element
$\sigma_\c\in\SL_2(\R)$ such that $\sigma_\c\infty\in\P^1(\R)$ is the
unique fixed point of~$\Gamma_\c$ in~$\P^1(\R)$ and such that
$$
\{\pm1\}
\sigma_\c^{-1}\Gamma_\c\sigma_\c=
\{\pm1\}
{\textstyle\bigl\{\smallmatrix 1b01\bigm|b\in\Z\bigr\}}.
$$
Such a $\sigma_\c$ exists and is unique up to multiplication
from the right by a matrix of the form~$\pm\smallmatrix1x01$ with
$x\in\R$; see Iwaniec~\citex{Iwaniec}{\S2.2}.  We define
$$
\eqalign{
q_\c\colon\hyp&\to\C\cr
z&\mapsto\exp(2\pi i\sigma_\c^{-1}z)}
$$
and
$$
\eqalign{
y_\c\colon\hyp&\to(0,\infty)\cr
z&\mapsto\Im\sigma_\c^{-1}z
=-{\log|q_\c(z)|\over2\pi}.}
$$

For all $\gamma\in\Gamma$, we write
$$
C_\c(\gamma)=|c|\quad\hbox{if }
\sigma_\c^{-1}\gamma\sigma_\c=\medmatrix abcd.
$$
Then we have
$$
\Gamma_\c=\{\gamma\in\Gamma\mid C_\c(\gamma)=0\}.
$$
It is known that the set $\{C_\c(\gamma)\mid\gamma\in\Gamma,
\gamma\not\in\Gamma_\c\}$ is bounded from below by a positive number,
and that if $\epsilon$ is a real number satisfying the inequality
$$
0<\epsilon\le\min_{\textstyle{\gamma\in\Gamma\atop
\gamma\not\in\Gamma_\c}}C_\c(\gamma),
\eqnumber{ineq:eps-sufficiently-small}
$$
then for all $z\in\hyp$ and~$\gamma\in\Gamma$ one has the implication
$$
y_\c(z)>1/\epsilon\hbox{ and }y_\c(\gamma z)>1/\epsilon
\implies\gamma\in\Gamma_\c.
$$
For any $\epsilon$ satisfying \eqref{ineq:eps-sufficiently-small}, the
image of the strip
$$
\{x+iy\mid 0\le x<1\hbox{ and }y>1/\epsilon\}\subset\hyp
$$
under the map
$$
\hyp\morphism{\sigma_\c}\hyp\longrightarrow\Gamma\backslash\hyp
$$
is an open disc $D_\c(\epsilon)$ around~$\c$, and the map $q_\c$
induces a chart on~$\Gamma\backslash\hyp$ identifying $D_\c(\epsilon)$
with the punctured disc $\{z\in\C\mid 0<|z|<\exp(-2\pi/\epsilon)\}$.
A compactification of~$\Gamma\backslash\hyp$ can be obtained by adding
a point for every cusp~$\c$ in such a way that $q_\c$ extends to a
chart with image equal to the disc
$\{z\in\C\mid|z|<\exp(-2\pi/\epsilon)\}$.  Let $\bar D_\c(\epsilon)$
denote the compactification of~$D_\c(\epsilon)$ obtained by adding the
boundary $\partial\bar D_\c(\epsilon)$ in~$\Gamma\backslash\hyp$ and
the cusp~$\c$.


\remark Let us fix a point~$w\in\hyp$ and write $\Gamma_w$ for the
stabilisator of~$w$ in~$\Gamma$.  The behaviour of $\gr_\Gamma(z,w)$
as $z\to w$ is
$$
\gr_\Gamma(z,w)={\#\Gamma_w\over2\pi}\log|z-w|\hbox{ as }z\to w.
$$
Furthermore, the behaviour of $\gr_\Gamma(z,w)$ as $z$ moves toward a
cusp~$\c$ of~$\Gamma$ is
$$
\gr_\Gamma(z,w)={1\over\vol_\Gamma}\log y_\c(z)+O(1)
\hbox{ as } y_\c(z)\to\infty.
$$

\subsection The Selberg--Harish-Chandra transform

\label{shc-transform}

Let $g\colon[1,\infty)\to\R$ be a smooth function with compact
support.  The {\it invariant integral operator attached to~$g$\/} is
the operator~$T_g$ defined on smooth functions $f\colon\hyp\to\C$ by
$$
(T_g f)(z)=\int_{w\in\hyp}g(u(z,w))f(w)\muhyp(w).
$$

The Laplace operator~$\Delta$ commutes with all such operators~$T_g$;
see Selberg~\citex{Selberg: Harmonic analysis}{pages 51--52} or
Iwaniec~\citex{Iwaniec}{Theorem~1.9}.  In fact, every eigenfunction
of~$\Delta$ is also an eigenfunction of all invariant integral
operators, and conversely; see Selberg~\citex{Selberg: Harmonic
analysis}{page~55} or Iwaniec~\citex{Iwaniec}{Theorems 1.14 and~1.15}.
The relation between the eigenvalues of~$\Delta$ and those of~$T_g$ is
given by the {\it Selberg--Harish-Chandra transform\/} of~$g$.  This
is a holomorphic function~$h$ defined by the following property.  Let
$f\colon\hyp\to\C$ be an eigenfunction of~$-\Delta$ with eigenvalue
$\lambda=s(1-s)$.  Then $f$ is also an eigenfunction of~$T_g$, and the
eigenvalue depends only on~$\lambda$; we can therefore define $h(s)$
uniquely such that
$$
-\Delta f=s(1-s)f\implies T_g f=h(s)f.
\eqnumber{selberg-transform}
$$
In particular, taking $f=1$, we see that
$$
h(0)=h(1)=2\pi\int_1^\infty g(u)du.
\eqnumber{eq:h01}
$$

The Selberg--Harish-Chandra transform can be identified with the
classical {\it Mehler--Fock transform\/}, defined as follows (see
Iwaniec~\citex{Iwaniec}{equation~$1.62'$}):
$$
h(s)=2\pi\int_1^\infty g(u)P_{s-1}(u)du.
\eqnumber{h-Legendre}
$$
Here $P_\nu$ is the Legendre function of the first kind of
degree~$\nu$; see Iwaniec~\citex{Iwaniec}{equation~1.43} or any book
on special functions, such as Erd\'elyi
et~al.~\citex{Erdelyi}{\S3.6.1}.  The function~$g$ can be recovered
from~$h$, and we call $g$ the {\it inverse Selberg--Harish-Chandra
transform\/} of~$h$.

The identity~\eqref{selberg-transform} holds more generally than just
for smooth functions~$g$ with compact support; see
Selberg~\citex{Selberg: Harmonic analysis}{pages 60--61}.  It will be
enough for us to state a slightly weaker, but more convenient
sufficient condition (cf.\ Selberg~\citex{Selberg: Harmonic
analysis}{page~72} or Iwaniec~\citex{Iwaniec}{equation 1.63}).  Let
$\epsilon>0$ and $\beta>1$, and let $h$ be a holomorphic function on
the strip $\{s\in\C\bigm|-\epsilon<\Re s<1+\epsilon\}$ such that
$h(s)=h(1-s)$ and such that $s\mapsto|h(s)||s(1-s)|^\beta$ is bounded
on this strip.  Then the inverse Selberg--Harish-Chandra transform~$g$
of~$h$ exists, and \eqref{selberg-transform} holds for the
pair~$(g,h)$.

\subsection Spectral theory of the Laplace operator for Fuchsian groups

Let $\Gamma$ be a cofinite Fuchsian group.  The spectrum
of~$-\Delta_\Gamma$ on~$\Ltwo(\Gamma\backslash\hyp)$ consists of a
discrete part and a continuous part.

The discrete spectrum consists of eigenvalues of~$-\Delta_\Gamma$ and
is of the form $\{\lambda_j\}_{j=0}^\infty$ with
$$
0=\lambda_0<\lambda_1\le\lambda_2\le\ldots,\quad
\lambda_j\to\infty\hbox{ as }j\to\infty.
$$
Let $\{\phi_j\}_{j=0}^\infty$ be a corresponding set of
eigenfunctions; these are called {\it automorphic forms of Maa\ss\ (of
weight~0)\/}.  We may and do assume that they are orthonormal with
respect to the inner product on~$\Ltwo(\Gamma\backslash\hyp)$.  For
each $j\ge0$, we define $s_j\in\C$ by
$$
\lambda_j=s_j(1-s_j),
$$
with $s_j\in[1/2,1]$ if $\lambda_j\le1/4$.  For $\lambda_j>1/4$, the
$s_j$ are only determined up to $s_j\leftrightarrow 1-s_j$.

The continuous part of the spectrum of~$-\Delta_\Gamma$ is the
interval $[1/4,\infty)$ with multiplicity equal to the number of cusps
of~$\Gamma$.  In particular, the continuous spectrum is absent if
$\Gamma$ has no cusps.  The continuous spectrum does not consist of
eigenvalues, but corresponds to ``wave packets'' that can be
constructed from {\it non-holomorphic Eisenstein series\/} or {\it
Eisenstein--Maa\ss\ series\/}, introduced by Maa\ss\ in~\cite{Maass}.
These series are defined as follows: for every cusp~$\c$
of~$\Gamma$ the series
$$
E_\c(z,s)=\sum_{\gamma\in\Gamma_\c\backslash\Gamma}
(y_\c(\gamma z))^s\quad
(z\in\hyp,s\in\C\hbox{ with }\Re s>1)
$$
converges uniformly on sets of the form $K\times\{s\in\C\mid\Re
s\ge\delta\}$ with $K$ a compact subset of~$\hyp$ and $\delta>1$.  In
particular, $E_\c(z,s)$ is a holomorphic function of~$s$.

A crucial ingredient in the spectral theory of automorphic forms is
the {\it meromorphic continuation of Eisenstein series\/}, due to
Selberg~\cite{Selberg: ICM}.  For proofs of this meromorphic
continuation and of the other properties of the Eisenstein series that
we use, we refer to Hejhal~\citex{Hejhal-II}{Chapter~VI, \S11}.
Different constructions of the meromorphic continuation can be found
in Faddeev~\citex{Faddeev}{\S4}, Hejhal~\citex{Hejhal-II}{Appendix~F}
or Iwaniec~\citex{Iwaniec}{Chapter~6}.

The meromorphic continuation of the Eisenstein series takes the
following form.  The functions $E_\c(z,s)$ can be continued to
functions of the form $E_\c(z,s)=H(z,s)/G(s)$, where $H$ is a smooth
function on~$\Gamma\backslash\hyp\times\C$ and both $G$ and~$H$ are
entire functions of~$s$.  These meromorphic continuations have a
finite number of simple poles on the segment~$(1/2,1]$ and no other
poles in $\{s\in\C\mid\Re s\ge1/2\}$, and satisfy a functional
equation, which we will not write down.  For all $s\in\C$ that is not
a pole, the function $z\mapsto E_\c(z,s)$ satisfies the differential
equation
$$
-\Delta_\Gamma E_\c(\blank,s)=s(1-s)E_\c(\blank,s).
$$
For $s\in\C$ with $\Re s=1/2$, the Eisenstein--Maa\ss\
series~$E_\c(\blank,s)$ are integrable, but not square-integrable, as
functions on~$\Gamma\backslash\hyp$.  In contrast, the ``wave
packets'' mentioned above are square-integrable.  They figure in
Theorem~\ref{theorem:spectral-representation} below, which is a
fundamental result in the theory of automorphic forms.


In the following theorem, and in the rest of the article, we will
consider integrals over the line $\Re s=1/2$.  For this we need an
orientation on this line; we fix one by requiring that the map
$t\mapsto 1/2+it$ from $\R$ with the usual orientation to the line
$\Re s=1/2$ preserves orientations.

\proclaimx Theorem (see Iwaniec~\citex{Iwaniec}{Theorems 4.7 and 7.3};
cf.\ Faddeev~\citex{Faddeev}{Theorem~4.1}). Every smooth and bounded
$\Gamma$-invariant function $f\colon\hyp\to\C$ has the spectral
representation
$$
f(z)=\sum_{j=0}^\infty b_j\phi_j(z)
+\sum_\c{1\over 4\pi i}\int_{\Re s=1/2} b_\c(s)\,
E_\c(z,s)ds,
\eqnumber{eq:spectral-representation}
$$
where $\c$ runs over the cusps of\/~$\Gamma$ and the
coefficients $b_j$ and $b_\c(s)$ are given by
$$
b_j=\int_{\Gamma\backslash\hyp}f\bar\phi_j\muhyp
\quad\hbox{and}\quad
b_\c(s)=\int_{\Gamma\backslash\hyp}f
\bar E_\c(\blank, s)\muhyp.
$$
The right-hand side of~\eqref{eq:spectral-representation} converges
to~$f$ in the Hilbert space $\Ltwo(\Gamma\backslash\hyp)$.  If in
addition the smooth $\Gamma$-invariant function $\Delta
f\colon\hyp\to\C$ is bounded, the convergence is uniform on compact
subsets of~$\hyp$.

\label{theorem:spectral-representation}

With regard to the spectral representations provided by this theorem,
the effect of the operator~$G_\Gamma$ from the introduction is as
follows: if $f$ has the spectral
representation \eqref{eq:spectral-representation}, then $G_\Gamma f$
has the corresponding spectral representation
$$
G_\Gamma f(z)=-\sum_{j=1}^\infty{b_j\over\lambda_j}\phi_j(z)
-\sum_\c{1\over 4\pi i}\int_{\Re s=1/2}{b_\c(s)\over s(1-s)}
E_\c(z,s)ds.
\eqnumber{eq:GGamma-spectral}
$$
(Note the absence of the eigenvalue $\lambda_0=0$.)

\def\topology{({\rm L}^\infty_{\rm loc},{\rm L}^2\cap {\rm
L}^\infty_{\rm loc})}

There is an analogous result
(Theorem~\ref{spectral-representation-kernel} below) for functions
on~$\hyp\times\hyp$ that are of the form
$\sum_{\gamma\in\Gamma}g(u(z,\gamma w))$, where
$g\colon[1,\infty)\to\R$ is the inverse Selberg--Harish-Chandra
transform of a function~$h$ as at the end of~\S\ref{shc-transform}.
The result involves a type of convergence that we now explain.  Let
$A$ be a filtered set, and let $\{K_a\}_{a\in A}$ be a family of
continuous functions
on~$\Gamma\backslash\hyp\times\Gamma\backslash\hyp$ that are
square-integrable in the second variable.  If $K$ is a function such
that for all compact subsets~$C$ of~$\Gamma\backslash\hyp$ we have
$$
\lim_{a\in A}\biggl(\sup_{z,w\in C}|K_a(z,w)-K(z,w)|
+\sup_{z\in C}\int_{w\in\Gamma\backslash\hyp}
|K_a(z,w)-K(z,w)|^2\muhyp(w)\biggr)=0,
$$
we say that the family of functions~$\{K_a\}_{a\in A}$ {\it converges
to~$K$ in the $\topology$-topology\/}.  In other words, this condition
means that the family converges uniformly on compact subsets
of~$\hyp\times\hyp$, and also with respect to the $\Ltwo$-norm in the
variable $w$, uniformly for $z$ in compact subsets
of~$\Gamma\backslash\hyp$.

\proclaimx Theorem (see Iwaniec~\citex{Iwaniec}{Theorem~7.4}).
Let $g\colon[1,\infty)\to\R$ be the inverse Selberg--Harish-Chandra
transform of a function~$h$ as at the end of~\S\ref{shc-transform}.
Then the function
$$
\eqalign{
K_g\colon\hyp\times\hyp&\longrightarrow\R\cr
(z,w)&\longmapsto\sum_{\gamma\in\Gamma}g(u(z,\gamma w))}
$$
is $\Gamma$-invariant with respect to both variables and admits the
spectral representation
$$
K_g(z,w)=\sum_{j=0}^\infty h(s_j)\phi_j(z)\bar\phi_j(w)
+\sum_\c{1\over 4\pi i}\int_{\Re s=1/2}
h(s)E_\c(z,s)\bar E_\c(w,s)ds,
\eqnumber{spectral-representation-kernel}
$$
where the expression on the right-hand side converges to~$K_g$ in the
$\topology$-topology.


\label{spectral-representation-kernel}

\subsection A point counting function

We fix a real number $U\ge1$, and define
$$
\eqalign{
g_U\colon[1,\infty)&\longrightarrow\R\cr
u&\longmapsto\cases{1& if $u\le U$;\cr 0& if $u>U$.}}
\eqnumber{definition-g_U}
$$
From \eqref{h-Legendre} and the formula for $\int_1^z P_\nu(w)dw$
found in Erd\'elyi et~al.~\citex{Erdelyi}{\S3.6.1, equation~8}, we see
that the Selberg--Harish-Chandra transform of~$g_U$ is
$$
h_U(s)=2\pi\sqrt{U^2-1}\,P^{-1}_{s-1}(U).
\eqnumber{eq:h-P}
$$
Here $P^\mu_\nu$ is the Legendre function of the first kind of
degree~$\nu$ and order~$\mu$; see \citex{Erdelyi}{\S3.2}.

Now let $\Gamma$ be a cofinite Fuchsian group.  We introduce the
following point counting function.  For any two points $z$, $w$
in~$\hyp$ and any $U\ge1$, we denote by $N_\Gamma(z,w,U)$ the number
of translates of~$w$ by elements of~$\Gamma$ lying in a disc
around~$z$ of radius~$r$ given by $\cosh(r)=U$, i.e.
$$
\eqalign{
N_\Gamma(z,w,U)&=\#\{\gamma\in\Gamma\mid u(z,\gamma w)\le U\}\cr
&=\sum_{\gamma\in\Gamma}g_U(u(z,w)).}
\eqnumber{eq:def-NGamma}
$$
This is $\Gamma$-invariant in $z$ and~$w$ separately.

\proclaim Lemma. Let $U\in[1,3]$, and let $s\in\C$ be such that
$s(1-s)(U-1)\in[0,1/2]$.  Then $h_U(s)$ is a real number satisfying
$$
(4\pi-8)(U-1)\le h_U(s)\le 8(U-1).
$$

\label{lemma:h_U-bound}

\proof This follows from \eqref{eq:h-P} and
Lemma~\ref{lemma:P-1-bound}.\endproof

\subsection Bounds on eigenfunctions

The convergence of the spectral
representation~\eqref{spectral-representation-kernel} can be deduced
from suitable bounds on the function
$$
\eqalign{
\Phi_\Gamma\colon\hyp\times[0,\infty)&\longrightarrow[0,\infty)\cr
(z,\lambda)&\longmapsto\sum_{j\colon\,\lambda_j\le\lambda}|\phi_j(z)|^2
+\sum_\c{1\over 4\pi i}
\int_{\textstyle{\Re s=1/2\atop s(1-s)\le\lambda}}
\bigl|E_\c(z,s)\bigr|^2ds.}
\eqnumber{eq:PhiGamma}
$$
We will prove a bound on~$\Phi_\Gamma$ which holds uniformly for all
subgroups~$\Gamma$ of finite index in a given Fuchsian
group~$\Gamma_0$.  This will give a similar uniformity in
Section~\ref{sec:explicit-bounds}.

\proclaim Lemma. Let $\Gamma$ be a cofinite Fuchsian group.  Then the
function $\Phi_\Gamma(z,\lambda)$ satisfies
$$
\Phi_\Gamma(z,\lambda)\le{\pi\over(2\pi-4)^2}N_\Gamma(z,z,17)\lambda
\quad\hbox{for all $z\in\hyp$ and all $\lambda\ge1/4$.}
$$

\label{eigenfunction-bound}

\proof Let $z\in\hyp$ and $\lambda\ge1/4$.  We put
$$
U=1+{1\over2\lambda}\in(1,3].
$$
From Bessel's inequality one can deduce (see
Iwaniec~\citex{Iwaniec}{\S7.2}) that
$$
\sum_{j\colon\,\lambda_j\le\lambda}|h_U(s_j)\phi_j(z)|^2
+\sum_\c{1\over 4\pi i}
\int_{\textstyle{\Re s=1/2\atop s(1-s)\le\lambda}}
\bigl|h_U(s)E_\c(z,s)\bigr|^2ds
\le\int_{w\in\Gamma\backslash\hyp}N_\Gamma(z,w,U)^2\muhyp(w).
$$
Using the definition~\eqref{eq:PhiGamma} of~$\Phi_\Gamma$ and the
bound $h_U(s)\ge(2\pi-4)/\lambda$ given by
Lemma~\ref{lemma:h_U-bound}, we deduce
$$
\Phi_\Gamma(z,\lambda)\le{\lambda^2\over(2\pi-4)^2}
\int_{w\in\Gamma\backslash\hyp}N_\Gamma(z,w,U)^2\muhyp(w).
$$
We rewrite the integral on the right-hand side by partial
``unfolding'' as follows (cf.\ Iwaniec~\citex{Iwaniec}{page~109}):
$$
\eqalign{
\int_{w\in\Gamma\backslash\hyp}N_\Gamma(z,w,U)^2\muhyp(w)
&=\sum_{\gamma,\gamma'\in\Gamma}\int_{w\in\Gamma\backslash\hyp}
g_U(z,\gamma' w)g_U(\gamma z,\gamma' w)\muhyp(w)\cr
&=\sum_{\gamma\in\Gamma}\int_{w\in\hyp}
g_U(z,w)g_U(\gamma z,w)\muhyp(w).}
$$
The last integral can be interpreted as the area of the intersection
of the discs of radius~$r$ around the points $z$ and~$\gamma z$
of~$\hyp$, where $\cosh r=U$.  By the triangle inequality for the
hyperbolic distance, this intersection is empty unless
$$
u(z,\gamma z)\le\cosh(2r)=2U^2-1;
$$
furthermore, the area of this intersection is at most
$2\pi(U-1)=\pi/\lambda$.  From this we deduce that
$$
\eqalign{
\int_{w\in\Gamma\backslash\hyp}N_\Gamma(z,w,U)^2\muhyp(w)
&\le{\pi\over\lambda}N_\Gamma(z,z,2U^2-1)}
$$
Since $2U^2-1\le17$, this proves the lemma.\endproof




\subsection The hyperbolic lattice point problem

\label{lattice-point-problem}

Let $\Gamma$ be a cofinite Fuchsian group.  The hyperbolic lattice
point problem for~$\Gamma$ is the following question: what is the
asymptotic behaviour of the point counting function $N_\Gamma(z,w,U)$
from~\eqref{eq:def-NGamma} as $U\to\infty$?  In contrast to the
Euclidean analogue of this question, about the number of points
in~$\Z^2$ lying inside a given circle in~$\R^2$, no elementary method
is known to even give the dominant term.  The difficulty is that for
circles in the hyperbolic plane of radius tending to infinity, the
circumference grows as fast as the enclosed area: the circumference of
a circle of radius~$r$ equals $2\pi\sinh(r)$, and the area of a disc
of radius~$r$ equals $2\pi(\cosh(r)-1)$.  In spite of this difficulty,
good estimates for~$N_\Gamma(z,w,U)$ can still be found, namely using
spectral theory on~$\Gamma\backslash\hyp$.

The strategy is to take suitable functions
$$
g_U^+, g_U^-\colon[1,\infty)\to\R
$$
with compact support, and to define functions $K_U^+$ and $K_U^-$
on~$\hyp\times\hyp$, invariant with respect to the action
of~$\Gamma$ on each of the two variables, by
$$
K_U^\pm(z,w)=\sum_{\gamma\in\Gamma}g_U^\pm(u(z,\gamma w)).
$$
This sum is finite because the functions~$g_U^\pm$ have compact
support.  We take the functions~$g_U^\pm$ such that for all
$z,w\in\hyp$ and $U>1$, we have the inequality
$$
K_U^-(z,w)\le N_\Gamma(z,w,U)\le K_U^+(z,w).
\eqnumber{inequality-N}
$$
Provided the Selberg--Harish-Chandra transforms $h_U^\pm$ of~$g_U^\pm$
satisfy the conditions of
Theorem~\ref{spectral-representation-kernel}, the functions $K_U^\pm$
have spectral representations
$$
K_U^\pm(z,w)=\sum_{j=0}^\infty h_U^\pm(s_j)\phi_j(z)\bar\phi_j(w)
+\sum_\c{1\over 4\pi i}\int_{\Re s=1/2}
h_U^\pm(s)E_\c(z,s)\bar E_\c(w,s)ds.
\eqnumber{KUpm-spectral}
$$
These spectral representations can then be used to find the asymptotic
behaviour of~$N_\Gamma(z,w,U)$ as~$U\to\infty$.

A reasonable choice at first sight would be to take for both $g_U^+$
and~$g_U^-$ the function~$g_U$ defined by~\eqref{definition-g_U}, so
that the inequalities in~\eqref{inequality-N} become equalities.
Unfortunately, the Selberg--Harish-Chandra transform~$h_U$ of~$g_U$
does not decay quickly enough as $|\Im s|\to\infty$ to give a spectral
representation of~$N_\Gamma(z,w,U)$ as in
Theorem~\ref{spectral-representation-kernel}.  We take instead
$$
g_U^+(u)=\cases{1& if $1\le u\le U$,\cr
{V-u\over V-U}& if $U\le u\le V$,\cr
0& if $V\le u$}
$$
and
$$
g_U^-(u)=\cases{1& if $1\le u\le T$,\cr
{U-u\over U-T}& if $T\le u\le U$,\cr
0& if $U\le u$}
$$
for certain $T$, $V$, depending on~$U$, with $1\le T<U<V$; see
Iwaniec~\citex{Iwaniec}{Chapter~12}.  Using \eqref{h-Legendre}, we
obtain
$$
\eqalign{
h_U^+(s)&=2\pi\int_1^V P_{s-1}(u){V-u \over V-U}du
-2\pi\int_1^V P_{s-1}(u){U-u\over V-U}du,\cr
h_U^-(s)&=2\pi\int_1^U P_{s-1}(u){U-u \over U-T}du
-2\pi\int_1^T P_{s-1}(u){T-u\over U-T}du.}
$$
Integrating by parts and applying the integral relation between the
Legendre functions $P_\nu$ and~$P^{-2}_\nu$ given in Erd\'elyi
et~al.~\citex{Erdelyi}{\S3.6.1, equation~8}, we get
$$
\eqalign{
h_U^+(s)&=2\pi{(V^2-1)P^{-2}_{s-1}(V)-(U^2-1)P^{-2}_{s-1}(U)
\over V-U},\cr
h_U^-(s)&=2\pi{(U^2-1)P^{-2}_{s-1}(U)-(T^2-1)P^{-2}_{s-1}(T)
\over U-T}.}
\eqnumber{eq:hUpm-Legendre}
$$
The dominant term in~\eqref{KUpm-spectral} as $U\to\infty$ comes from
the eigenvalue $\lambda_0=0$, corresponding to $s_0=1$.  It follows
from~\eqref{eq:h01} or the formula $P^{-2}_{0}(u)=(u-1)/(2u+2)$ that
$$
h_U^+(1)=2\pi(U-1)+\pi(V-U)
\quad\hbox{and}\quad
h_U^-(1)=2\pi(U-1)-\pi(U-T).
\eqnumber{h_U^+(1)}
$$

Let a real number $\delta\ge1$ be given.  We fix parameters
$\alpha^+$, $\alpha^-$, $\beta^+$ and~$\beta^-$ satisfying
$$
\alpha^\pm\in(0,1/2),\quad
\beta^\pm>0,\quad
\beta^-\le{\delta^{1+\alpha^-}\over\delta+1}.
\eqnumber{eq:beta-crit}
$$
We choose $T$ and~$V$ as functions of~$U$ as follows:
$$
T(U)=U-\beta^- U^{-1-\alpha^-}(U^2-1),\qquad
V(U)=U+\beta^+ U^{-1-\alpha^+}(U^2-1).
\eqnumber{choice-T-V}
$$
The last inequality in~\eqref{eq:beta-crit} ensures that if
$U\ge\delta$, then $T(U)\ge1$.


For later use, we will keep the parameters $\alpha^\pm$
and~$\beta^\pm$ variable for greater flexibility.  To obtain the best
known error bound in the hyperbolic lattice point problem, the right
choice is $\alpha^\pm=1/3$, so that
$$
V-U\sim\beta^+ U^{2/3}
\quad\hbox{and}\quad
U-T\sim\beta^- U^{2/3}
\quad\hbox{as }U\to\infty.
$$
This choice leads to the following theorem.

\proclaimx Theorem (Huber~\citex{Huber}{Satz~B},
Patterson~\citex{Patterson}{Theorem~2}, Selberg; see
Iwaniec~\citex{Iwaniec}{Theorem~12.1}).  Let $\Gamma$ be a cofinite
Fuchsian group.  For all $z,w\in\hyp$, the point counting function
$N_\Gamma$ satisfies
$$
N_\Gamma(z,w,U)=\sum_{j\colon\,2/3<s_j\le 1}2^{s_j}\sqrt{\pi}
{\Gamma\bigl(s_j-{1\over2}\bigr)\over\Gamma(s_j+1)}
\phi_j(z)\bar\phi_j(w)U^{s_j}+O(U^{2/3})\quad\hbox{as }U\to\infty,
$$
with an implied constant depending on~$\Gamma$ and the points $z$
and~$w$.

\label{lattice-point-estimate}

In particular, since $|\phi_0|^2$ is the constant
function~$1/{\vol_\Gamma}$, this shows that
$$
N_\Gamma(z,w,U)\sim{2\pi(U-1)\over\vol_\Gamma}\quad\hbox{as }U\to\infty.
$$
Since $2\pi(U-1)$ is the area of a disc of radius~$r$ with $\cosh
r=U$, Theorem~\ref{lattice-point-estimate} implies that this area is
asymptotically equivalent to the number of lattice points inside the
disc times the area of a fundamental domain for the action
of~$\Gamma$, which is the intuitively expected result.

\section An approximate spectral representation of the Green function

\label{section:approximate-spectral-representation}

Let $\Gamma$ be a cofinite Fuchsian group.  The Green function
of~$\Gamma$ {\it formally\/} has the spectral representation
$$
\gr_\Gamma(z,w)\buildrel?\over=
-\sum_{j=1}^\infty{1\over\lambda_j}\phi_j(z)\bar\phi_j(w)
-\sum_\c{1\over 4\pi i}\int_{\Re s=1/2}{1\over s(1-s)}
E_\c(z,s)\bar E_\c(w,s)ds.
$$
The problem is that this expansion does not converge.  Neither should
one be tempted to write the Green function by ``averaging'' $\gr_\hyp$
as a (likewise divergent) sum
$$
\gr_\Gamma(z,w)\buildrel?\over=
\sum_{\gamma\in\Gamma}\gr_\hyp(z,\gamma w).
$$
However, both of these divergent expressions have at least some value
as guiding ideas for what follows.  In fact, we will bound
$\gr_\Gamma(z,w)$ by means of certain functions
$R_{\Gamma,\delta}^\pm(z,w)$, defined in~\eqref{definition-R} below,
that reflect the above formal spectral representation of~$\gr_\Gamma$.

\subsection A construction of the Green function using the resolvent
kernel

\label{construction-green-function}

We will give a construction of the Green function of~$\Gamma$ using
the family of auxiliary functions
$$
\eqalign{
g_a\colon(1,\infty)&\longrightarrow[0,\infty)\cr
u&\longmapsto{1\over2\pi}Q_{a-1}(u)}
$$
for $a\ge1$, where $Q_\nu$ is the Legendre function of the second kind
of degree~$\nu$; see Erd\'elyi et~al.~\citex{Erdelyi}{\S3.6.1}.
By \citex{Erdelyi}{\S3.6.2, equation~20}, we have
$$
Q_0(u)={1\over2}\log{u+1\over u-1},
$$
which shows that $g_1$ equals the function~$L$
from~\eqref{definition-L}.  By \eqref{h-Legendre} and
\citex{Erdelyi}{\S3.12, equation~4}, the Selberg--Harish-Chandra
transform of~$g_a$ is
$$
\eqalign{
h_a(s)&=\int_1^\infty P_{s-1}(u)Q_{a-1}(u)du\cr
&={1\over(a-s)(a-1+s)}\cr
&={1\over s(1-s)+a(a-1)}.}
$$

Given a real number $\sigma<1/2$, we consider the strip
$$
S_\sigma=\{s\in\C\mid\sigma\le\Re s\le 1-\sigma\}.
\eqnumber{eq:strip}
$$

\proclaim Lemma. (a)\enspace For all $a,b>1$ and all
$\sigma\in(1-\min\{a,b\},1/2)$, the function
$$
s\mapsto\bigl|h_a(s)-h_b(s)\bigr|\bigl|s(1-s)\bigr|^2
$$
is bounded on~$S_\sigma$.
\smallskip\noindent
(b)\enspace Let $\sigma\in(0,1/2)$.  There exist real numbers
$(C_{a,b,\sigma})_{a,b>1}$, with $C_{a,b,\sigma}\to0$ as both $a$
and~$b$ tend to~$1$, such that
$$
\bigl|h_a(s)-h_b(s)\bigr|\le C_{a,b,\sigma}\bigl|s(1-s)\bigr|^{-2}
\quad\hbox{for all }s\in S_\sigma.
$$

\label{lemma:prop-ha-hb}

\proof Both claims are easily deduced from the expression
$$
h_a(s)-h_b(s)={b(b-1)-a(a-1)\over(a-s)(a-1+s)(b-s)(b+1-s)}.
$$
Details are left to the reader.\endproof

For all $a>1$, the sum $\sum_{\gamma\in\Gamma}g_a(u(z,\gamma w))$
converges uniformly on compact subsets of~$\hyp\times\hyp$ not
containing any points of the form $(z,\gamma z)$ and defines a
continuous function that is square-integrable in each variable; see
Fay~\citex{Fay}{Theorem~1.5}.
We can therefore define
$$
\eqalign{
K_a^\Gamma\colon\{(z,w)\in\hyp\times\hyp\mid z\not\in\Gamma w\}
&\longrightarrow\R\cr
(z,w)&\longmapsto\sum_{\gamma\in\Gamma}g_a(u(z,\gamma w))-c_a,}
\eqnumber{definition-KGamma}
$$
where
$$
c_a={2\pi\over \vol_\Gamma}\int_1^\infty g_a(u)du
={1\over\vol_\Gamma}h_a(1)={1\over\vol_\Gamma a(a-1)}.
$$
The constant~$c_a$ is such that the integral of~$K_a^\Gamma$
over~$\Gamma\backslash\hyp$ with respect to each of the variables
vanishes.  Up to this constant, $K_a^\Gamma$ is the resolvent kernel
with parameter~$a$.

It is known that the resolvent kernel admits a meromorphic
continuation in the variable~$a$.  The following result can be
interpreted as the statement that $-\gr_\Gamma$ is the constant term
in the Laurent expansion of the resolvent kernel at $a=1$.

\proclaim Proposition.  The family of functions
$\{-K^\Gamma_a\}_{a>1}$ converges to the Green function~$\gr_\Gamma$
in the $\topology$-topology.

\label{prop:gr-Gamma}

\proof
It follows from Lemma~\ref{lemma:prop-ha-hb}(a) that for all $a,b>1$,
the function $g_a-g_b$ satisfies the conditions of
Theorem~\ref{spectral-representation-kernel}.  The function
$$
(K^\Gamma_a-K^\Gamma_b)(z,w)=\sum_{\gamma\in\Gamma}
\bigl(g_a(u(z,\gamma w))-g_b(u(z,\gamma w))\bigr)-c_a+c_b
$$
therefore has the spectral representation
$$
\eqalign{
(K^\Gamma_a-K^\Gamma_b)(z,w)&=\sum_{j=1}^\infty(h_a(s_j)-h_b(s_j))
\phi_j(z)\bar\phi_j(w)\cr
&\qquad+\sum_\c{1\over 4\pi i}\int_{\Re s=1/2}(h_a(s)-h_b(s))
E_\c(z,s)\bar E_\c(w,s)ds,}
\eqnumber{eq:Ka-Kb-spectral}
$$
where the right-hand side converges to~$K^\Gamma_a-K^\Gamma_b$ in the
$\topology$-topology.  (Note that the eigenvalue $\lambda_0=0$ has
disappeared because of the definition of~$c_a$.)  In particular,
$K^\Gamma_a-K^\Gamma_b$ extends to a continuous function
on~$\hyp\times\hyp$ that is $\Gamma$-invariant with respect to both
variables.

We claim that $\{K^\Gamma_a-K^\Gamma_b\}_{a,b>1}$ converges to~0 in
the $\topology$-topology as $a,b\searrow1$.  In particular, this
implies that $\{K^\Gamma_a\}_{a>1}$ converges to a symmetric
continuous function
on~$\Gamma\backslash\hyp\times\Gamma\backslash\hyp$ that is
square-integrable with respect to each variable separately.  We fix
$\sigma\in(0,1/2)$ be such that the spectrum of~$-\Delta_\Gamma$ is
contained in $\{0\}\cup[\sigma(1-\sigma),\infty)$.

First we show that $\{K^\Gamma_a-K^\Gamma_b\}_{a,b>1}$ converges to
zero uniformly on compact subsets of~$\hyp\times\hyp$.
Lemma~\ref{lemma:prop-ha-hb}(b) implies
$$
\eqalign{
\bigl|K^\Gamma_a-K^\Gamma_b\bigr|(z,w)&\le
\sum_{j=1}^\infty|h_a(s_j)-h_b(s_j)|\cdot|\phi_j(z)\bar\phi_j(w)|\cr
&\qquad+\sum_\c{1\over 4\pi i}\int_{\Re s=1/2}|h_a(s)-h_b(s)|\cdot
\bigl|E_\c(z,s\bigr)\bar E_\c(w,s)\bigr|ds\cr
&\le C_{a,b,\sigma}\Biggl(\sum_{j=1}^\infty
\bigl(s_j(1-s_j)\bigr)^{-2}|\phi_j(z)\bar\phi_j(w)|\cr
&\hskip18mm+\sum_\c{1\over 4\pi i}\int_{\Re s=1/2}
\bigl(s(1-s)\bigr)^{-2}\bigl|E_\c(z,s)\bar E_\c(w,s)\bigr|ds\Biggr).}
$$
By the Cauchy--Schwarz inequality and Lemma~\ref{eigenfunction-bound},
the right-hand side converges to~0 uniformly on compact subsets
of~$\hyp\times\hyp$, as claimed.

Next we show that $\{K^\Gamma_a-K^\Gamma_b\}_{a,b>1}$ converges to
zero with respect to the $\Ltwo$-norm in the variable~$w$, uniformly
for $z$ in compact subsets of~$\hyp$.  From \eqref{eq:Ka-Kb-spectral},
Plancherel's theorem and Lemma~\ref{lemma:prop-ha-hb}(b), we deduce
$$
\eqalign{
\int_{w\in\Gamma\backslash\hyp}
\bigl|K^\Gamma_a-K^\Gamma_b\bigr|^2(z,w)\muhyp(w)
&=\sum_{j=1}^\infty\bigl|h_a(s_j)-h_b(s_j)\bigr|^2|\phi_j(z)|^2\cr
&\qquad+\sum_\c{1\over 4\pi i}\int_{\Re s=1/2}
\bigl|h_a(s)-h_b(s)\bigr|^2\bigl|E_\c(z,s)\bigr|^2ds\cr
&\le C_{a,b,\sigma}^2\Biggl(\sum_{j=1}^\infty
\bigl(s_j(1-s_j)\bigr)^{-4}|\phi_j(z)|^2\cr
&\hskip18mm+\sum_\c{1\over 4\pi i}
\int_{\Re s=1/2}\bigl(s(1-s)\bigr)^{-4}
\bigl|E_\c(z,s)\bigr|^2 ds\Biggr)}
$$
for real numbers $C_{a,b,\sigma}$ with $C_{a,b,\sigma}\to0$ as
$a,b\searrow1$.  By Lemma~\ref{eigenfunction-bound}, the other factor
on the right-hand side is bounded on compact subsets of~$\hyp$.  This
implies that the right-hand side converges to~0 uniformly on compact
subsets of~$\hyp$ as $a,b\searrow1$, as claimed.

The defining property~\eqref{selberg-transform} of the
Selberg--Harish-Chandra transform implies that if $f$ is a smooth,
bounded, $\Gamma$-invariant function on~$\hyp$, with spectral
representation \eqref{eq:spectral-representation}, then
$$
\int_{w\in\Gamma\backslash\hyp}K^\Gamma_a(z,w)f(w)\muhyp(w)
=\sum_{j=1}^\infty b_j h_a(s_j)\phi_j(z)
+\sum_\c{1\over 4\pi i}\int_{\Re s=1/2} b_\c(s)
h_a(s)E_\c(z,s)ds.
$$
Taking the limit, using the $\Ltwo$-convergence that we just proved
and applying \eqref{eq:GGamma-spectral}, we get
$$
\eqalign{
\int_{w\in\Gamma\backslash\hyp}
\lim_{a\searrow1}K^\Gamma_a(z,w)f(w)\muhyp(w)&=
\lim_{a\searrow1}\int_{w\in\Gamma\backslash\hyp}
K^\Gamma_a(z,w)f(w)\muhyp(w)\cr
&=\sum_{j=1}^\infty {b_j\over s_j(1-s_j)}\phi_j(z)
+\sum_\c{1\over 4\pi i}\int_{\Re s=1/2}
{b_\c(s)\over s(1-s)}E_\c(z,s)ds\cr
&=-G_\Gamma f(z)\cr
&=-\int_{w\in\Gamma\backslash\hyp}\gr_\Gamma(z,w)f(w)\muhyp(w).}
$$
Since the set of smooth and bounded functions is dense
in~$\Ltwo(\Gamma\backslash\hyp)$, this proves that the limit of the
convergent family of functions $\{K^\Gamma_a\}_{a>1}$ equals
$-\gr_\Gamma$.\endproof

\subsection Proof of the approximate spectral representation

\label{subsec:approximate-spectral-representation}

We now exploit the estimates for the hyperbolic lattice point problem
given in~\S\ref{lattice-point-problem}.  We choose parameters
$\alpha^\pm$ and~$\beta^\pm$ satisfying~\eqref{eq:beta-crit}.  Using
these, we define functions $T(U)$, $V(U)$, $g_U^\pm(u)$, $h_U^\pm(s)$
and $K_U^\pm(z,w)$ as in~\S\ref{lattice-point-problem}.  Furthermore,
we fix a real number $\delta>1$.  We write $\Gamma$ as the disjoint
union of subsets $\Pi_{\Gamma,\delta}(z,w)$
and~$\Lambda_{\Gamma,\delta}(z,w)$ defined by
$$
\eqalign{
\Pi_{\Gamma,\delta}(z,w)&=\{\gamma\in\Gamma\mid u(z,\gamma w)\le\delta\},\cr
\Lambda_{\Gamma,\delta}(z,w)&=\{\gamma\in\Gamma\mid u(z,\gamma w)>\delta\}.}
$$
We define
$$
I_\delta^\pm(s)={1\over2\pi}\int_\delta^\infty{h_U^\pm(s)\over U^2-1}dU
\quad\hbox{for }0<\Re s<1,
\eqnumber{definition-I}
$$
$$
R_{\Gamma,\delta}^\pm(z,w)=
\sum_{j=1}^\infty I_\delta^\pm(s_j)\phi_j(z)\bar\phi_j(w)
+\sum_\c{1\over 4\pi i}\int_{\Re s=1/2} I_\delta^\pm(s)
E_\c(z,s)\bar E_\c(w,s)ds,
\eqnumber{definition-R}
$$
$$
q_{\Gamma,\delta}^+={1\over\vol_\Gamma}
\biggl({\beta^+\over2\alpha^+\delta^{\alpha^+}}
-\log{\delta+1\over 2}\biggr),
\quad
q_{\Gamma,\delta}^-=-{1\over\vol_\Gamma}
\biggl({\beta^-\over2\alpha^-\delta^{\alpha^-}}
+\log{\delta+1\over 2}\biggr).
\eqnumber{eq:definition-q}
$$

The intuition behind the following theorem is that although the Green
function $\gr_\Gamma$ does not admit a spectral representation, it can
be bounded (after removing the logarithmic singularity) by functions
that do admit spectral representations.  The terms
$q_{\Gamma,\delta}^\pm$ below correspond to the eigenvalue~$0$, while
the terms $R_{\Gamma,\delta}^\pm(z,w)$ correspond to the non-zero part
of the spectrum.

\proclaim Theorem.  Let $\Gamma$ be a cofinite Fuchsian group.  For
all $\delta>1$ and for every choice of the parameters $\alpha^\pm$
and~$\beta^\pm$ satisfying~\eqref{eq:beta-crit}, the Green function
of\/~$\Gamma$ satisfies the inequalities
$$
-q_{\Gamma,\delta}^+ - R_{\Gamma,\delta}^+(z,w)
\le\gr_\Gamma(z,w)+\sum_{\gamma\in\Pi_{\Gamma,\delta}(z,w)}
\bigl(L(u(z,\gamma w))-L(\delta)\bigr)
\le -q_{\Gamma,\delta}^- - R_{\Gamma,\delta}^-(z,w).
$$

\label{big-theorem}

\proof For any $U\ge\delta$, the inequality~\eqref{inequality-N}
implies that the number of elements
$\gamma\in\Lambda_{\Gamma,\delta}(z,w)$ with $u(z,\gamma w)\le U$ can
be bounded as
$$
A(U)\le\#\{\gamma\in\Lambda_{\Gamma,\delta}(z,w)\mid
u(z,\gamma w)\le U\}\le B(U),
\eqnumber{ineq:AUB}
$$
where the functions $A,B\colon [\delta,\infty)\to\R$ are defined by
$$
A(U)=K_U^-(z,w)-\#\Pi_{\Gamma,\delta}(z,w)
\quad\hbox{and}\quad
B(U)=K_U^+(z,w)-\#\Pi_{\Gamma,\delta}(z,w).
$$
The functions $A$ and~$B$ are continuous and increasing.  The
estimates from~\S\ref{lattice-point-problem} imply that they are
bounded linearly in~$U$ as $U\to\infty$, with an implied
constant depending on the group~$\Gamma$, the points $z$ and~$w$ and
the functions $T$ and~$V$.

Let $\{h_a\}_{a>1}$, $\{g_a\}_{a>1}$ and $\{K^\Gamma_a\}_{a>1}$ be as
in~\S\ref{construction-green-function}.  For all $a>1$, applying
partial summation and \eqref{ineq:AUB} gives
$$
-\int_\delta^\infty g_a'(U)A(U)dU\le
\sum_{\gamma\in\Lambda_{\Gamma,\delta}(z,w)} g_a(u(z,\gamma w))
\le-\int_\delta^\infty g_a'(U)B(U)dU.
$$
Using the definition~\eqref{definition-KGamma} of~$K^\Gamma_a$, we
deduce the upper bound
$$
K^\Gamma_a(z,w)\le\sum_{\gamma\in\Pi_{\Gamma,\delta}(z,w)}g_a(u(z,\gamma w))
-\int_\delta^\infty g_a'(U)B(U)dU
-{2\pi\over \vol_\Gamma}\int_1^\infty g_a(u)du.
$$
The definition of~$B$ implies
$$
\eqalign{
\int_\delta^\infty g_a'(U)B(U)dU
&=\int_\delta^\infty g_a'(U)K_U^+(z,w)dU
-\#\Pi_{\Gamma,\delta}(z,w)\int_\delta^\infty g_a'(U)dU\cr
&=\int_\delta^\infty g_a'(U)\Bigl(K_U^+(z,w)
-{2\pi\over\vol_\Gamma}(U-1)\Bigr)dU
+{2\pi\over\vol_\Gamma}\int_\delta^\infty g_a'(U)(U-1)dU\cr
&\qquad+\#\Pi_{\Gamma,\delta}(z,w) g_a(\delta).}
$$
Using integration by parts, we rewrite the second integral in the last
expression as follows:
$$
\eqalign{
\int_\delta^\infty g_a'(U)(U-1)dU
&=\int_1^\infty g_a'(U)(U-1)dU
-\int_1^\delta g_a'(U)(U-1)dU\cr
&=-\int_1^\infty g_a(U)dU-\int_1^\delta g_a'(U)(U-1)dU.}
$$
We can now rewrite our upper bound for~$K^\Gamma_a(z,w)$ as
$$
\eqalign{
K^\Gamma_a(z,w)&\le\sum_{\gamma\in\Pi_{\Gamma,\delta}(z,w)}
\bigl(g_a(u(z,w))-g_a(\delta)\bigr)-\int_\delta^\infty g_a'(U)
\Bigl(K_U^+(z,w)-{2\pi\over\vol_\Gamma}(U-1)\Bigr)dU\cr
&\qquad+{2\pi\over\vol_\Gamma}\int_1^\delta g_a'(U)(U-1)dU.}
$$
Lemma~\ref{lemma:Q'-bound} implies
$$
{1\over2\pi}\biggl({2\over u+1}\biggr)^{a-1}{1\over u^2-1}
\le g_a'(u)\le 0,
$$
and equality holds for $a=1$.  By the dominated convergence theorem,
we may take the limit $a\searrow1$ inside the integrals.
Together with Proposition~\ref{prop:gr-Gamma}, this leads to
$$
\eqalign{
\gr_\Gamma(z,w)+\sum_{\gamma\in\Pi_{\Gamma,\delta}(z,w)}
\bigl(L(u(z,\gamma w))-L(\delta)\bigr)
&\ge-{1\over2\pi}\int_\delta^\infty
\Bigl(K_U^+(z,w)-{2\pi\over \vol_\Gamma}(U-1)\Bigr)
{dU\over U^2-1}\cr
&\qquad+{1\over \vol_\Gamma}\log{\delta+1\over2}.}
$$
In the integral, we insert the spectral
representation~\eqref{KUpm-spectral} of~$K_U^+$, the
formula~\eqref{h_U^+(1)} for $h_U^+(1)$ and the fact that
$|\phi_0|^2=1/{\vol_\Gamma}$.  We then interchange the resulting sums
and integrals with the integral over~$U$; this is permitted because
the double sums and integrals converge absolutely, as one deduces from
Lemma~\ref{eigenfunction-bound} and
Theorem~\ref{lattice-point-estimate}.  This yields
$$
{1\over2\pi}\int_\delta^\infty
\Bigl(K_U^+(z,w)-{2\pi\over \vol_\Gamma}(U-1)\Bigr)
{dU\over U^2-1}=R_{\Gamma,\delta}^+(z,w)
+{1\over 2\vol_\Gamma}\int_\delta^\infty{V-U\over U^2-1}dU.
$$
Finally, we note that
$$
\eqalign{
\int_\delta^\infty{V-U\over U^2-1}dU
&=\beta^+\int_\delta^\infty U^{-1-\alpha^+} dU\cr
&={\beta^+\over\alpha^+\delta^{\alpha^+}}.}
$$
This proves the lower bound of the theorem.  The proof of the upper
bound is similar.\endproof

\remark The only inequality responsible for the fact that the
inequalities in Theorem~\ref{big-theorem} are not equalities
is~\eqref{ineq:AUB}.

\section Explicit bounds

\label{sec:explicit-bounds}

\subsection Bounds on $h_U^\pm(s)$ and $I_\delta^\pm(s)$

We keep the notation
of~\S\ref{subsec:approximate-spectral-representation}.  In addition,
we choose real numbers $\sigma^\pm$ such that
$$
0<\alpha^+<\sigma^+<1/2
\quad\hbox{and}\quad
0<\alpha^-<\sigma^-<1/2.
$$

Let $s$ be in the strip~$S_{\sigma^+}$ defined by~\eqref{eq:strip},
and let $p_{\sigma^+}(u)$ be the elementary function defined
by~\eqref{eq:psig} below.  From \eqref{eq:hUpm-Legendre},
Corollary~\ref{cor:bound-P2-regular} and~\eqref{choice-T-V}, we obtain
$$
\eqalign{
|h_U^+(s)|&\le2\pi{(V^2-1)\bigl|P^{-2}_{s-1}(V)\bigr|
+(U^2-1)\bigl|P^{-2}_{s-1}(U)\bigr|\over V-U}\cr
&\le 2\pi\bigl|s(1-s)\bigr|^{-5/4}
{p_{\sigma^+}(V)+p_{\sigma^+}(U)\over V-U}\cr
&=2\pi\bigl|s(1-s)\bigr|^{-5/4}
{\bigl(p_{\sigma^+}(V)+p_{\sigma^+}(U)\bigr)U^{1+\alpha^+}
\over\beta^+(U^2-1)}.}
$$
Similarly, for $s\in S_{\sigma^-}$,
$$
|h_U^-(s)|\le 2\pi\bigl|s(1-s)\bigr|^{-5/4}
{\bigl(p_{\sigma^-}(U)+p_{\sigma^-}(T)\bigr)U^{1+\alpha^-}
\over\beta^-(U^2-1)}.
$$
Substituting this in the definition~\eqref{definition-I} of~$I$, we
obtain
$$
|I_\delta^+(s)|\le D^+_\delta\bigl|s(1-s)\bigr|^{-5/4}
\quad\hbox{and}\quad
|I_\delta^-(s)|\le D^-_\delta\bigl|s(1-s)\bigr|^{-5/4},
\eqnumber{eq:Ipmregular}
$$
where
$$
\eqalign{
D^+_\delta&={1\over\beta^+}
\int_\delta^\infty{\bigl(p_{\sigma^+}(V)+p_{\sigma^+}(U)\bigr)
U^{1+\alpha^+}\over(U^2-1)^2}dU,\cr
D^-_\delta&={1\over\beta^-}
\int_\delta^\infty{\bigl(p_{\sigma^-}(U)+p_{\sigma^-}(T)\bigr)
U^{1+\alpha^-}\over(U^2-1)^2}dU.}
\eqnumber{eq:Dpm}
$$


\subsection Bounds on~$\gr_\Gamma$

\proclaim Theorem.  Let $\Gamma$ be a Fuchsian group.  Let $\delta>1$
and~$\eta\in(0,1/4]$ be real numbers such that the spectrum
of~$-\Delta_\Gamma$ is contained in $\{0\}\cup[\eta,\infty)$.  Let
$\sigma^+$, $\sigma^-$, $\alpha^+$, $\alpha^-$, $\beta^+$, $\beta^-$
be real numbers satisfying \eqref{eq:beta-crit} and the inequalities
$$
0<\alpha^+<\sigma^+<1/2,\quad
0<\alpha^-<\sigma^-<1/2
\quad\hbox{and}\quad
\sigma^\pm(1-\sigma^\pm)\le\eta.
$$
Then the Green function of\/~$\Gamma$ satisfies the inequalities
$$
A\le\gr_\Gamma(z,w)+\sum_{\gamma\in\Pi_{\Gamma,\delta}(z,w)}
\bigl(L(u(z,\gamma w))-L(\delta)\bigr)\le B
\quad\hbox{for all }z,w\in\hyp,
$$
where
$$
\eqalign{
A&=-q_{\Gamma,\delta}^+ - D^+_\delta{\pi\over(2\pi-4)^2}
\biggl({\eta^{-5/4}\over4}+4\sqrt2\biggr)
{N_\Gamma(z,z,17)+N_\Gamma(w,w,17)\over2},\cr
B&=-q_{\Gamma,\delta}^- + D^-_\delta{\pi\over(2\pi-4)^2}
\biggl({\eta^{-5/4}\over4}+4\sqrt2\biggr)
{N_\Gamma(z,z,17)+N_\Gamma(w,w,17)\over2}.}
$$

\label{theorem:explicit}

\proof In view of Theorem~\ref{big-theorem}, we have to bound the
absolute values of the functions~$R_{\Gamma,\delta}^\pm(z,w)$
from~\eqref{definition-R}.  Applying the triangle inequality and the
Cauchy--Schwarz inequality, we see that
$$
\bigl|R_{\Gamma,\delta}^\pm(z,w)\bigr|\le{S^\pm(z)+S^\pm(w)\over2},
$$
where $S^+$ and~$S^-$ are defined by
$$
S^\pm(z)=\sum_{j=1}^\infty|I_\delta^\pm(s_j)||\phi_j(z)|^2
+\sum_\c{1\over 4\pi i}\int_{\Re s=1/2}|I_\delta^\pm(s)|
\bigl|E_\c(z,s)\bigr|^2ds.
$$
Let $\Phi_\Gamma(z,\lambda)$ be as in~\eqref{eq:PhiGamma}.  Applying
\eqref{eq:Ipmregular}, we obtain (with
$\partial\Phi_\Gamma/\partial\lambda$ taken in a distributional sense)
$$
\eqalign{
S^\pm(z)/D^\pm_\delta
&\le\sum_{j=1}^\infty\lambda_j^{-5/4}|\phi_j(z)|^2
+\sum_\c{1\over 4\pi i}\int_{\Re s=1/2}\bigl(s(1-s)\bigr)^{-5/4}
\bigl|E_\c(z,s)\bigr|^2ds\cr
&\le\eta^{-5/4}\sum_{j\colon\,\lambda_j\le1/4}|\phi_j(z)|^2
+\int_{1/4}^\infty\lambda^{-5/4}
{\partial\Phi_\Gamma\over\partial\lambda}(z,\lambda)d\lambda\cr
&=\eta^{-5/4}\Phi_\Gamma(z,1/4)
+\left[\lambda^{-5/4}\Phi_\Gamma(z,\lambda)\right]_{\lambda=1/4}^\infty
+{5\over4}\int_{1/4}^\infty\lambda^{-9/4}\Phi_\Gamma(z,\lambda)d\lambda\cr
&=\bigl(\eta^{-5/4}-2^{5/2}\bigr)\Phi_\Gamma(z,1/4)
+{5\over4}\int_{1/4}^\infty\lambda^{-9/4}\Phi_\Gamma(z,\lambda)d\lambda.}
$$
The bound on~$\Phi_\Gamma(z,\lambda)$ given by
Lemma~\ref{eigenfunction-bound} implies
$$
\eqalign{
S^\pm(z)&\le
D^\pm_\delta{\pi\over(2\pi-4)^2}N_\Gamma(z,z,17)
\biggl(\bigl(\eta^{-5/4}-2^{5/2}\bigr)\cdot{1\over4}
+{5\over4}\int_{1/4}^\infty\lambda^{-9/4}\lambda\,d\lambda\biggr)\cr
&=D^\pm_\delta{\pi\over(2\pi-4)^2}N_\Gamma(z,z,17)
\biggl({\eta^{-5/4}\over4}+4\sqrt2\biggr).}
$$
This proves the theorem.\endproof



\medbreak\noindent {\it Proof of
  Theorem~\ref{theorem:introduction}.\/}\enspace Let the notation be
as in the theorem; we may assume $\eta\le1/4$.  We apply
Theorem~\ref{theorem:explicit} to~$\Gamma$, with parameters
$\sigma^\pm$, $\alpha^\pm$ and~$\beta^\pm$ depending only on~$\eta$
and not on~$\Gamma$.  It is clear that the factor $1/{\vol_\Gamma}$
occurring in the definition~\eqref{eq:definition-q} is bounded
by~$1/{\vol_{\Gamma_0}}$, and that $N_\Gamma(z,z,17)$ is bounded
by~$N_{\Gamma_0}(z,z,17)$.  It remains to remark that
$N_{\Gamma_0}(z,z,17)$ is bounded on~$Y_0$.\endproof

The bounds given by Theorem~\ref{theorem:explicit} are easy to make
explicit.  First, the real numbers $D^\pm_\delta$ from~\eqref{eq:Dpm}
can be bounded in elementary ways or approximated by numerical
integration.  Second, a straightforward computation shows that for
$z=x+iy\in\hyp$ and $\gamma=\smallmatrix abcd\in\SL_2(\R)$, we have
$$
u(z,\gamma z)={1\over2}\biggl(
(a-cx)^2+\left({b+(a-d)x-cx^2\over y}\right)^2
+(cy)^2+(d+cx)^2\biggr).
\eqnumber{eq:uabcd}
$$
This can be used for concrete groups~$\Gamma$ to find an upper bound
on $N_\Gamma(z,z,U)$ for $U>1$.

\subsection Example: congruence subgroups of $\SL_2(\Z)$

\label{subsection:example}

Let us consider the case where $\Gamma_0=\SL_2(\Z)$.  We will make the
bounds from Theorem~\ref{theorem:introduction} explicit for congruence
subgroups $\Gamma\subseteq\SL_2(\Z)$.  By our convention for
integration on $\Gamma\backslash\hyp$ if $-1\in\Gamma$, we have
$$
{1\over\vol_\Gamma}\le{1\over\vol_{\SL_2(\Z)}}={\pi\over6}.
$$
We choose
$$
\delta=2.
$$
Selberg conjectured in~\cite{Selberg: Fourier coefficients} that the
least non-zero eigenvalue~$\lambda_1$ of~$-\Delta_\Gamma$ is at
least~1/4, and he proved that $\lambda_1\ge3/16$.  The sharpest result
known so far, due to Kim and Sarnak \citex{Kim}{Appendix~2}, is that
$\lambda_1\ge(25/64)(1-25/64)=975/4096$.  We may therefore take
$$
\eta=975/4096.
$$

We now consider the point counting function~$N_{\SL_2(\Z)}(z,z,U)$
defined by~\eqref{eq:def-NGamma} on a rectangle of the form
$$
R=\{x+iy\in\hyp\mid
x_{\min}\le x\le x_{\max}, y_{\min}\le y\le y_{\max}\}
$$
for given real numbers $x_{\min}<x_{\max}$ and $0<y_{\min}<y_{\max}$.
The function $z\mapsto N_{\SL_2(\Z)}(z,z,U)$ on~$R$ is clearly bounded
from above by the number of matrices $\gamma=\smallmatrix
abcd\in\SL_2(\Z)$ such that for some $z\in R$, the inequality
$$
u(z,\gamma z)\le U
\eqnumber{ineq:u}
$$
holds.  We now show how to enumerate these matrices.  We distinguish
the cases $c=0$ and $c\ne0$.  We assume (after multiplying by $-1$ if
necessary) that $a=d=1$ in the first case, and that $c>0$ in the
second case.  The total number of matrices~$\gamma$ as above is then
twice the number produced by our enumeration.

In the case $c=0$, by \eqref{eq:uabcd}, the inequality~\eqref{ineq:u}
reduces to
$$
1+\textfrac1/2(b/y)^2\le U.
$$
This implies
$$
|b|\le y_{\max}\sqrt{2U-2}.
$$
In the case $c>0$, it follows from \eqref{eq:uabcd} and~\eqref{ineq:u}
that
$$
\displaylines{
|c|\le\sqrt{2U}/y_{\min},\cr
-\sqrt{2U}+cx_{\min}\le a\le\sqrt{2U}+cx_{\max},\cr
-\sqrt{2U}-cx_{\max}\le d\le\sqrt{2U}-cx_{\min}.}
$$
Since $c\ne0$, the coefficients $a$, $c$, $d$ and the condition
$ad-bc=1$ determine~$b$.  If $\gamma$ is a matrix obtained in this
way, we compute the minimum of~$u(z,\gamma z)$ for $z\in R$
using~\eqref{eq:uabcd} to decide whether there exists a point $z\in R$
satisfying~\eqref{ineq:u}.

Let $Y_0$ denote the compact subset in $\SL_2(\Z)\backslash\hyp$ which
is the image of the rectangle
$$
\{x+iy\in\hyp\mid-1/2\le x\le1/2
\hbox{ and }\sqrt{3}/2\le y\le 2\}.
$$
This is the complement of a disc around the unique cusp
of~$\SL_2(\Z)$.  Dividing this rectangle into $100\times100$ small
rectangles and bounding $N_{\SL_2(\Z)}(z,z,U)$ on each of them as
described above, we get
$$
N_{\Gamma_0}(z,z,17)\le216
\quad\hbox{for all }z\in Y_0.
$$
Given this upper bound for $N_{\SL_2(\Z)}(z,z,17)$, some
experimentation leads to the following values for the parameters:
$$
\displaylines{
\alpha^+=0.0366,\quad\beta^+=2.72,\quad\sigma^+=0.306,\cr
\alpha^-=2.96\cdot10^{-3},\quad\beta^-=0.668,\quad\sigma^-=0.250.}
$$
With these choices, a numerical calculation gives
$$
q_{\Gamma,\delta}^+<69.0,\quad
q_{\Gamma,\delta}^->-216,\quad
D^+_\delta<18.5,\quad
D^-_\delta<9.61.
$$
This implies the following explicit bounds on Green functions of
congruence subgroups.

\proclaim Theorem. Let $\Gamma$ be a congruence subgroup
of~$\SL_2(\Z)$.  Then for all $z,w\in\hyp$ whose images
in~$\SL_2(\Z)\backslash\hyp$ lie in~$Y_0$, we have
$$
-2.87\cdot10^4\le
\gr_\Gamma(z,w)+\sum_{\textstyle{\gamma\in\Gamma\atop u(z,\gamma w)\le2}}
\biggl(L(u(z,\gamma w))-{\log3\over4\pi}\biggr)
\le 1.51\cdot10^4.
$$

\section Extension to neighbourhoods of the cusps

\label{section:extension-cusps}

The bounds given by Theorem~\ref{theorem:explicit} do not have the
right asymptotic behaviour when $z$ or~$w$ are near the cusps
of~$\Gamma$.  This means that we have to do some more work to find
suitable bounds on the Green function $\gr_\Gamma(z,w)$ in this case.


Let $D$ and~$\bar D$ denote the open and closed unit discs in~$\C$,
respectively.  We recall that the Poisson kernel on~$D$ is defined by
$$
\eqalign{
P(\zeta)&={1-|\zeta|^2\over|1-\zeta|^2}\cr
&=1+\sum_{n=1}^\infty \zeta^n+\sum_{n=1}^\infty\bar\zeta^n.}
$$
We will use the notation
$$
\tilde P(t,\zeta)=P(\exp(2\pi i t)\zeta).
$$

\proclaim Lemma. The Poisson kernel satisfies
$$
\int_0^1 \tilde P(a,\zeta) \tilde P(-a,\eta)da
=P(\zeta\eta)\quad\hbox{for all }\zeta,\eta\in D
\eqnumber{eq:P-convolution}
$$
and
$$
\tilde P(t,\zeta)={d\over dt}\left(t+
{1\over2\pi i}\left(
\log{1-\exp(-2\pi i t)\bar\zeta\over 1-\exp(2\pi i t)\zeta}
-\log{1-\bar\zeta\over 1-\zeta}
\right)\right)\quad\hbox{for all }\zeta\in D.
\eqnumber{eq:P-deriv}
$$

\proof The first claim can be verified in several ways, for example
using the residue theorem, Fourier series, or the fact that the
Poisson kernel solves the Laplace equation with Dirichlet boundary
conditions.  The second claim is straightforward to check.\endproof

Let $\gr_{\bar D}$ denote the Green function for the Laplace operator
on~$\bar D$; this is an integral kernel for the Poisson equation
$\Delta f=g$ with boundary condition $f=0$ on $\partial\bar D$.  It is
given explicitly by
$$
\gr_{\bar D}(\zeta,\eta)=
{1\over2\pi}\log\left|{\zeta-\eta\over 1-\zeta\bar\eta}\right|
\quad\hbox{for all $\zeta,\eta\in D$ with $\zeta\ne\eta$}.
$$


For all $\xi\in D$ and $t\in\R$, we write
$$
\lambda(\xi,t)={1\over 2\pi i}
\bigl(\log(1-\exp(-2\pi i t)\xi)-\log(1-\exp(2\pi i t)\xi)\bigr).
$$

\proclaim Lemma. The function~$\lambda(\xi,t)$ satisfies
$$
\biggl|\int_0^t{\lambda(\xi,y)\over y}dy
+{1\over2}\log(1-\xi)\biggr|\le{1\over 12t}
\quad\hbox{for all }t>0.
$$

\label{lemma:lambda-integral}

\proof We expand $\lambda(\xi,t)$ for $\xi\in D$ in a Fourier series:
$$
\eqalign{
\lambda(\xi,t)&={1\over2\pi i}\left(
\sum_{n=1}^\infty {\xi^n\exp(2\pi i n t)\over n}
-\sum_{n=1}^\infty {\xi^n\exp(-2\pi i n t)\over n}\right)\cr
&={1\over\pi}\sum_{n=1}^\infty{\xi^n\sin(2\pi n t)\over n}.}
$$
This implies
$$
\eqalign{
\int_0^t{\lambda(\xi,y)\over y}dy&=
{1\over\pi}\sum_{n=1}^\infty{\xi^n\over n}
\int_0^t{\sin(2\pi n y)\over y}dy\cr
&={1\over\pi}\sum_{n=1}^\infty{\xi^n\over n}
\int_0^{2\pi n t}{\sin x\over x}dx\cr
&={1\over\pi}\sum_{n=1}^\infty{\xi^n\over n}
(\si(0)-\si(2\pi n t)).}
$$
Here $\si(y)$ is the sine integral function normalised such that
$\lim_{y\to\infty}\si(y)=0$:
$$
\si(y)=\int_y^\infty{\sin x\over x}dx.
$$
It is known that
$$
\si(0)={\pi\over2}
\quad\hbox{and}\quad
|{\si(x)}|\le{1\over x}\quad\hbox{for all }x>0.
$$
From this we get
$$
\eqalign{
\int_0^ t{\lambda(\xi,y)\over y}dy&=
{1\over2}\sum_{n=1}^\infty{\xi^n\over n}
-{1\over\pi}\sum_{n=1}^\infty{\xi^n\si(2\pi n t)\over n}\cr
&=-{1\over2}\log(1-\xi)
-{1\over\pi}\sum_{n=1}^\infty{\xi^n\si(2\pi n t)\over n}}
$$
and
$$
\eqalign{
\left|\sum_{n=1}^\infty{\xi^n\si(2\pi n t)\over n}\right|
&\le\sum_{n=1}^\infty{1\over n(2\pi n t)}\cr
&={1\over 2\pi t}\sum_{n=1}^\infty{1\over n^2}\cr
&={1\over 2\pi t}{\pi^2\over6}.}
$$
This proves the claim.\endproof

For $\delta>1$ and $u>1$, we write
$$
J_\delta(u)=\max\{0,L(u)-L(\delta)\}.
$$
For $\xi\in D$, $\delta>1$ and $\epsilon>0$, we write
$$
N_{\delta,\epsilon}(\xi)=\int_{t\in\R}J_\delta
\biggl(1+{(\epsilon t)^2\over2}\biggr)P(\exp(2\pi i t)\xi)dt.
$$

\proclaim Lemma. The function $N_{\delta,\epsilon}$ satisfies
$$
\left|N_{\delta,\epsilon}(\xi)-{1\over\epsilon}\cdot
{2\over\pi}\arctan\sqrt{\delta-1\over2}
+{1\over2\pi}\log|1-\xi|\right|
\le\epsilon r_\delta
\quad\hbox{for all }\xi\in D,
$$
where
$$
r_\delta={1\over 24\pi}\biggl(\sqrt{2\over\delta-1}
+\arctan\sqrt{\delta-1\over2}\biggr).
\eqnumber{eq:definition-Rdelta}
$$

\label{lemma:bounds-N}

\proof We note that
$$
1+{(\epsilon t)^2\over 2}\le\delta
\iff
|t|\le\tau,
$$
where
$$
\tau={\sqrt{2\delta-2}\over\epsilon}.
$$
By the definition of~$J_\delta$, this gives
$$
N_{\delta,\epsilon}(\xi)=\int_{-\tau}^\tau\biggl(
L\biggl(1+{(\epsilon t)^2\over2}\biggr)
-L(\delta)\biggr)\tilde P(t,\xi)dt.
$$
Using \eqref{eq:P-deriv}, integrating by parts, and taking the
contributions for positive and negative~$t$ together, we obtain
$$
\eqalign{
N_{\delta,\epsilon}(\xi)&=-\int_{-\tau}^\tau \epsilon^2 t L'(1+(\epsilon t)^2/2)
\biggl(t+{1\over2\pi i}\biggl(
\log{1-\exp(-2\pi i t)\bar\xi\over 1-\exp(2\pi i t)\xi}
-\log{1-\bar\xi\over 1-\xi}
\biggr)\biggr)dt\cr
&=-\int_0^\tau \epsilon^2 t L'(1+(\epsilon t)^2/2)
\biggl(2t+{1\over2\pi i}\biggl(
\log{1-\exp(-2\pi i t)\bar\xi\over 1-\exp(2\pi i t)\xi}
-\log{1-\exp(2\pi i t)\bar\xi\over 1-\exp(-2\pi i t)\xi}\biggr)\biggr)dt\cr
&=-\int_0^\tau \epsilon^2 t L'(1+(\epsilon t)^2/2)
\bigl(2t+\lambda(\xi,t)
+\lambda(\bar\xi,t)\bigr)dt.}
$$
Using the definition~\eqref{definition-L} of~$L$ and rearranging
gives
$$
\eqalign{
N_{\delta,\epsilon}(\xi)&={1\over2\pi}\int_0^\tau \epsilon^2 t
\biggl({1\over(\epsilon t)^2}-{1\over 4+(\epsilon t)^2}\biggr)
\bigl(2t+\lambda(\xi,t)+\lambda(\bar\xi,t)\bigr)dt\cr
&={1\over2\pi}\int_0^\tau
\biggl(2-{2(\epsilon t)^2\over 4+(\epsilon t)^2}\biggr)dt
+{1\over2\pi}\int_0^\tau
\biggl(1-{(\epsilon t)^2\over 4+(\epsilon t)^2}\biggr)
{\lambda(\xi,t)+\lambda(\bar\xi,t)\over t}dt\cr
&={1\over2\pi}\int_0^\tau{8\over 4+(\epsilon t)^2}dt
+{1\over2\pi}\int_0^\tau{4\over 4+(\epsilon t)^2}
{\lambda(\xi,t)+\lambda(\bar\xi,t)\over t}dt.}
\eqnumber{eq:N-final}
$$
We consider the two integrals in the last expression one by one.  As
for the first integral, we have
$$
\eqalign{
\int_0^\tau{8\over 4+(\epsilon t)^2}dt&=
{2\over\epsilon}\int_0^{\epsilon\tau/2}{2\over 1+x^2}dx\cr
&={4\over\epsilon}\arctan{\epsilon\tau\over2}\cr
&={4\over\epsilon}\arctan\sqrt{\delta-1\over2}.}
\eqnumber{eq:N-arctan}
$$
As for the second integral in~\eqref{eq:N-final}, let us write for
convenience
$$
I_\xi=\int_0^\tau{4\over 4+(\epsilon t)^2}
{\lambda(\xi,t)+\lambda(\bar\xi,t)\over t}dt
$$
and
$$
\Lambda_\xi(t)=\int_0^t{\lambda(\xi,y)+\lambda(\bar\xi,y)\over y}dy
+\log|1-\xi|.
$$
Then we have
$$
\Lambda_\xi'(t)={\lambda(\xi,t)+\lambda(\bar\xi,t)\over t}
\quad\hbox{and}\quad
\Lambda_\xi(0)=\log|1+\xi|.
$$
Integration by parts gives
$$
I_\xi=-\log|1-\xi|+{4\over 4+(\epsilon\tau)^2}\Lambda_\xi(\tau)
+\int_0^\tau{8\epsilon^2t\over(4+(\epsilon t)^2)^2}\Lambda_\xi(t)dt.
$$
By Lemma~\ref{lemma:lambda-integral}, it follows that
$$
\bigl|I_\xi+\log|1-\xi|\bigr|\le
+{4\over 4+(\epsilon\tau)^2}{1\over 6\tau}
+\int_0^\tau{8\epsilon^2t\over(4+(\epsilon t)^2)^2}{1\over 6t}dt.
$$
The integral can be evaluated by elementary means, and the result is
$$
\eqalign{
\bigl|I_\xi+\log|1-\xi|\bigr|&\le
{1\over6\tau}+{\epsilon\over12}\arctan{\epsilon\tau\over2}\cr
&={\epsilon\over12}\biggl(\sqrt{2\over\delta-1}
+\arctan\sqrt{\delta-1\over2}\biggr).}
$$
Combining this with \eqref{eq:N-final} and~\eqref{eq:N-arctan} proves
the claim.\endproof


\proclaim Lemma. {Let $\Gamma$ be a cofinite Fuchsian group, and let
  $\delta>1$ and $\epsilon'>\epsilon>0$ be real numbers satisfying the
  inequalities
$$
\bigl(\delta+\sqrt{\delta^2-1}\bigr)^{1/2}\epsilon'\le
\min_{\textstyle{\gamma\in\Gamma\atop\gamma\not\in\Gamma_\c}}
C_\c(\gamma)
\quad\hbox{and}\quad
\bigl(\delta+\sqrt{\delta^2-1}\bigr)\epsilon\le\epsilon'.
$$
(a)\enspace For all $z,w\in\hyp$ with $y_\c(z)\ge1/\epsilon'$ and
$y_\c(w)\ge1/\epsilon'$ and all $\gamma\in\Gamma$, we have
$$
u(z,\gamma w)<\delta\implies\gamma\in\Gamma_\c.
$$
(b)\enspace For all $z,w\in\hyp$ such that $y_\c(z)\ge1/\epsilon$ and
such that the image of~$w$ in~$\Gamma\backslash\hyp$ lies outside
$D_\c(\epsilon')$, and for all $\gamma\in\Gamma$, we have
$u(z,\gamma w)\ge\delta$.}

\label{lemma:bla}

\proof Let $z$, $w$ and~$\gamma$ be as in~(a).  We write
$$
\sigma_\c^{-1}\gamma\sigma_\c=\gamma'=\medmatrix abcd.
$$
Suppose $\gamma\not\in\Gamma_\c$.  Then our assumptions imply
$$
|c|/\epsilon'\ge\bigl(\delta+\sqrt{\delta^2-1}\bigr)^{1/2}.
\eqnumber{ineq:bla}
$$
We have
$$
\eqalign{
u(z,\gamma w)&=u(\sigma_\c^{-1}z,\sigma_\c^{-1}\gamma w)\cr
&=u(\sigma_\c^{-1}z,\gamma'\sigma_\c^{-1}w)\cr
&=1+{|\sigma_\c^{-1}z-\gamma'\sigma_\c^{-1}w|^2\over
2(\Im\sigma_\c^{-1}z)(\Im\gamma'\sigma_\c^{-1}w)}\cr
&\ge1+{(\Im\sigma_\c^{-1}z-\Im\gamma'\sigma_\c^{-1}w)^2\over
2(\Im\sigma_\c^{-1}z)(\Im\gamma'\sigma_\c^{-1}w)}\cr
&={1\over2}\biggl(
{\Im\sigma_\c^{-1} z\over\Im\gamma'\sigma_\c^{-1} w}
+{\Im\gamma'\sigma_\c^{-1} w\over\Im\sigma_\c^{-1} z}\biggr)\cr
&={1\over2}\biggl(
{y_\c(z)|c\sigma_\c^{-1}w+d|^2\over y_\c(w)}
+{y_\c(w)\over y_\c(z)|c\sigma_\c^{-1}w+d|^2}\biggr).}
$$
From~\eqref{ineq:bla}, we deduce
$$
\eqalign{
{y_\c(z)|c\sigma_\c^{-1}w+d|^2\over y_\c(w)}
&\ge{y_\c(z)(c y_\c(w))^2\over y_\c(w)}\cr
&=c^2 y_\c(z) y_\c(w)\cr
&\ge (|c|/\epsilon')^2\cr
&\ge\delta+\sqrt{\delta^2-1}.}
$$
Using the fact that the function $x\mapsto x+x^{-1}$ is increasing for
$x\ge1$, we obtain
$$
\eqalign{
u(z,\gamma w)
&\ge{1\over2}\biggl(\bigl(\delta+\sqrt{\delta^2-1}\bigr)
+{1\over\delta+\sqrt{\delta^2-1}}\biggr)\cr
&=\delta.}
$$
This proves (a).

Now let $z$, $w$ and~$\gamma$ be as in~(b).  Our assumption that the
image of~$w$ in~$\Gamma\backslash\hyp$ lies outside~$D_\c(\epsilon')$
implies
$$
y_\c(\gamma w)\le 1/\epsilon'
$$
and hence
$$
\eqalign{
{y_\c(z)\over y_\c(\gamma w)}&\ge{\epsilon'\over\epsilon}\cr
&\ge\delta+\sqrt{\delta^2-1}.}
$$
Using the fact that the function $x\mapsto x^{-1}$ is increasing for
$x\ge1$ as in the proof of~(a), we get
$$
\eqalign{
u(z,\gamma w)&=u(\sigma_\c^{-1}z,\sigma_\c^{-1}\gamma w)\cr
&\ge{1\over2}\biggl({y_\c(z)\over y_\c(\gamma w)}
+{y_\c(\gamma w)\over y_\c(z)}\biggr)\cr
&\ge{1\over2}\biggl(\bigl(\delta+\sqrt{\delta^2-1}\bigr)
+{1\over\delta+\sqrt{\delta^2-1}}\biggr)\cr
&=\delta.}
$$
This proves (b).\endproof


In the following proposition, we extend our bounds on~$\gr_\Gamma$ to
the neighbourhoods $D_\c(\epsilon_\c)$ of the cusps.  We make the
following abuse of notation: for $z\in\hyp$ and $S$ a subset
of~$\Gamma\backslash\hyp$, we write $z\in S$ if the image of~$z$
in~$\Gamma\backslash\hyp$ lies in~$S$.

\proclaim Proposition.  {Let $\Gamma$ be a cofinite Fuchsian group,
and let $\delta$ be a real number with $\delta>1$.  For every
cusp~$\c$ of~$\Gamma$, let $\epsilon_\c'>\epsilon_\c>0$ be real
numbers satisfying the inequalities
$$
\epsilon_\c'\bigl(\delta+\sqrt{\delta^2-1}\bigr)^{1/2}\le
\min_{\textstyle{\gamma\in\Gamma\atop\gamma\not\in\Gamma_\c}}
C_\c(\gamma),
\quad\hbox{and}\quad
\bigl(\delta+\sqrt{\delta^2-1}\bigr)\epsilon_\c\le\epsilon_\c'
$$
and small enough such that the discs $D_\c(\epsilon_\c')$ are pairwise
disjoint.  Let
$$
Y=(\Gamma\backslash\hyp)\mathbin{\big\backslash}
\bigsqcup_\c D_\c(\epsilon_\c).
$$
Let $A$ and~$B$ be real numbers satisfying
$$
A\le \gr_\Gamma(z,w)+
\sum_{\textstyle{\gamma\in\Gamma\atop u(z,\gamma w)\le\delta}}
\bigl(L(u(z,w))-L(\delta)\bigr)\le B
\quad\hbox{for all }z,w\in Y.
\eqnumber{ineq:assumption-AB}
$$
(a)\enspace If $\c$ is a cusp such that $z\in D_\c(\epsilon_\c)$,
$w\in Y$ and $w\not\in D_\c(\epsilon_\c')$, we have
$$
A\le\gr_\Gamma(z,w)-
{1\over\vol_\Gamma}\log(\epsilon_\c y_\c(z))\le B.
$$
(a$'$)\enspace If $\c$ is a cusp such that $w\in D_\c(\epsilon_\c)$,
$z\in Y$ and $z\not\in D_\c(\epsilon_\c')$, we have
$$
A\le\gr_\Gamma(z,w)-
{1\over\vol_\Gamma}\log(\epsilon_\c y_\c(w))\le B.
$$
(b) If $\c$, $\fd$ are two distinct cusps such that $z\in
D_\c(\epsilon_\c)$ and $w\in D_\fd(\epsilon_\fd)$, we have
$$
A\le\gr_\Gamma(z,w)-
{1\over\vol_\Gamma}\log(\epsilon_\c y_\c(z))-
{1\over\vol_\Gamma}\log(\epsilon_\c y_\c(w))\le B.
$$
(c) If $\c$ is a cusp such that $z,w\in D_\c(\epsilon_\c')$, we have
$$
\eqalign{
\tilde A_\c\le\gr_\Gamma(z,w)
&-\#(\Gamma\cap\{\pm1\})\cdot{1\over2\pi}\log|q_\c(z)-q_\c(w)|\cr
&-{1\over\vol_\Gamma}\log(\epsilon_\c' y_\c(z))
-{1\over\vol_\Gamma}\log(\epsilon_\c' y_\c(w))
\le\tilde B_\c,}
$$
where $\tilde A_\c$ and $\tilde B_\c$ are defined using the
function~$r_\delta$ from~\eqref{eq:definition-Rdelta} by
$$
\eqalign{
\tilde A_\c&=A
+\#(\Gamma\cap\{\pm1\})\biggl[{1\over\epsilon_\c'}
\biggl(1-{2\over\pi}\arctan\sqrt{\delta-1\over2}\biggr)
-\epsilon_\c' r_\delta\biggr],\cr
\tilde B_\c&=B
+\#(\Gamma\cap\{\pm1\})\biggl[{1\over\epsilon_\c'}
\biggl(1-{2\over\pi}\arctan\sqrt{\delta-1\over2}\biggr)
+\epsilon_\c' r_\delta\biggr].}
$$
}


\proof In view of \S\ref{cusps} (or Lemma~\ref{lemma:bla}(a)), the
discs $D_\c(\epsilon_\c')\supset D_\c(\epsilon_\c)$ are well defined.
Furthermore, the assumption that the discs $D_\c(\epsilon_\c')$ are
pairwise disjoint implies that for every cusp~$\c$, the boundaries of
$\bar D_\c(\epsilon_\c)$ and $\bar D_\c(\epsilon_\c')$ are contained
in~$Y$.

Let us prove part (a).  We keep $w\in Y$ fixed and consider
$\gr_\Gamma(z,w)$ as a function of~$z\in D_\c(\epsilon_\c)$.  The
defining properties of~$\gr_\Gamma$ imply
$$
\gr_\Gamma(z,w)={1\over\vol_\Gamma}
\log(\epsilon_\c y_\c(z))+h_w(z)
\quad\hbox{for all }z\in D_\c(\epsilon_\c),
$$
with $h_w$ a real-valued harmonic function on $\bar
D_\c(\epsilon_\c)$.  By construction, $h_w(z)$ coincides with
$\gr_\Gamma(z,w)$ for $z$ on the boundary of~$\bar D_\c(\epsilon_\c)$.
This implies
$$
h_w(z)=\int_0^1\gr_\Gamma(\sigma_\c(a+i/\epsilon_\c),w)
\tilde P(-a,q_\c(z)\exp(2\pi/\epsilon_\c))da.
$$
By Lemma~\ref{lemma:bla}(b), there are no $\gamma\in\Gamma$ such that
$u(z,\gamma w)<\delta$.  By the assumption~\eqref{ineq:assumption-AB}
on $A$ and~$B$, we conclude
$$
A\le h_w(z)\le B\quad\hbox{for all }z\in D_\c(\epsilon).
$$
This proves (a).  Part (a$'$) is equivalent to~(a) by symmetry, and
(b) is proved in a similar way.

It remains to prove part~(c).  We identify $\bar D_\c(\epsilon_\c')$ with
the closed unit disc~$\bar D$ via the map
$$
\eqalign{
\bar D_\c(\epsilon_\c')&\isom\bar D\cr
z&\longmapsto\zeta_z=q_\c(z)\exp(2\pi/\epsilon_\c').}
$$
Let $\gr_{\bar D_\c(\epsilon_\c')}$ be the Green function for the
Laplace operator on~$\bar D_\c(\epsilon_\c')$, given in terms of~$\gr_D$
by
$$
\gr_{\bar D_\c(\epsilon_\c')}(z,w)=
\#(\Gamma\cap\{\pm1\})\cdot
\gr_D(\zeta_z,\zeta_w).
$$
The factor $\#(\Gamma\cap\{\pm1\})$ arises because of how we defined
integration on~$\Gamma\backslash\hyp$ in
Section~\ref{sec:introduction}.

Fixing $w$ and considering $\gr_\Gamma$ as a function of~$z$, we have
$$
\gr_\Gamma(z,w)=\gr_{\bar D_\c(\epsilon_\c')}(z,w)
+{1\over\vol_\Gamma}\log(\epsilon_\c' y_\c(z))
+h_w(z)\quad\hbox{for all }z\in D_\c(\epsilon_\c'),
$$
where $h_w$ is a real-valued harmonic function on $\bar
D_\c(\epsilon_\c')$.  By construction, $h_w(z)$ coincides with
$\gr_\Gamma(z,w)$ for $z$ on the boundary of~$\bar
D_\c(\epsilon_\c')$.  This implies
$$
h_w(z)=\int_0^1 \gr_\Gamma(\sigma_\c(a+i/\epsilon_\c'),w)
\tilde P(-a,\zeta_z)da.
$$
Applying the same argument to $\gr_\Gamma(\sigma_\c(a+i/\epsilon_\c'),w)$
as a function of~$w$, we obtain
$$
\gr_\Gamma(z,w)=\gr_{\bar D_\c(\epsilon_\c')}(z,w)
+{1\over\vol_\Gamma}\log(\epsilon_\c' y_\c(z))
+{1\over\vol_\Gamma}\log(\epsilon_\c' y_\c(w))
+K(z,w),
\eqnumber{eq:grGamma-K}
$$
where $K$ is the function on $D_\c(\epsilon_\c')\times D_\c(\epsilon_\c')$
defined by
$$
K(z,w)=\int_{a=0}^1\int_{b=0}^1
\gr_\Gamma(\sigma_\c(a+i/\epsilon_\c'),\sigma_\c(b+i/\epsilon_\c'))
\tilde P(-b,\zeta_w)\tilde P(-a,\zeta_z)db\,da.
$$
In \eqref{ineq:assumption-AB}, we may replace $\Gamma$ by~$\Gamma_\c$
in view of Lemma~\ref{lemma:bla}(a), i.e., we have
$$
A\le
\gr_\Gamma(z,w)+\sum_{\textstyle{\gamma\in\Gamma_\c\atop u(z,\gamma w)\le\delta}}
\bigl(L(u(z,\gamma w))-L(\delta)\bigr)
\le B
\quad\hbox{for all }z,w\in \partial\bar D_\c(\epsilon_\c').
$$
Substituting this in the definition of~$K(z,w)$ and ``unfolding'' the
action of~$\Gamma_\c$, we get
$$
A\le K(z,w) + \#(\Gamma\cap\{\pm1\})M(\zeta_z,\zeta_w)\le B,
\eqnumber{ineq:A-M<=K<=B-M}
$$
where $M$ is the function on $D\times D$ defined by
$$
M(\zeta,\eta)=\int_{a=0}^1\int_{b\in\R}
J_\delta(u(\sigma_\c(a+i/\epsilon_\c'),\sigma_\c(b+i/\epsilon_\c')))
\tilde P(-b,\eta)\tilde P(-a,\zeta)db\,da.
$$
Making the change of variables
$$
b=a+t
$$
and noting that
$$
\eqalign{
u(\sigma_\c(a+i/\epsilon_\c'),\sigma_\c(b+i/\epsilon_\c'))
&=u(a+i/\epsilon_\c',b+i/\epsilon_\c')\cr
&=1+{(b-a)^2\over2/\epsilon_\c'^2}\cr
&=1+{(\epsilon_\c' t)^2\over2},}
$$
we obtain
$$
M(\zeta,\eta)=\int_{a=0}^1\int_{t\in\R}
J_\delta\biggl(1+{(\epsilon_\c' t)^2\over2}\biggr)
\tilde P(-a-t,\eta)\tilde P(-a,\zeta)dt\,da.
$$
Interchanging the order of integration, noting that
$$
\tilde P(-a-t,\eta)=\tilde P(a,\exp(2\pi i t)\bar\eta)
$$
and using \eqref{eq:P-convolution}, we simplify this to
$$
\eqalign{
M(\zeta,\eta)&=\int_{t\in\R}J_\delta\biggl(1+{(\epsilon_\c' t)^2\over2}\biggr)
\tilde P(t,\zeta\bar\eta)dt\cr
&=N_{\delta,\epsilon_\c'}(\zeta\bar\eta).}
$$
Applying Lemma~\ref{lemma:bounds-N}, we conclude
from~\eqref{ineq:A-M<=K<=B-M} that
$$
\eqalign{
K(z,w)&\le B+\#(\Gamma\cap\{\pm1\})\biggl(
{1\over2\pi}\log|1-\zeta_z\bar\zeta_w|
-{1\over\epsilon_\c'}\cdot{2\over\pi}\arctan\sqrt{\delta-1\over2}
+\epsilon_\c' r_\delta\biggr),\cr
K(z,w)&\ge A+\#(\Gamma\cap\{\pm1\})\biggl(
{1\over2\pi}\log|1-\zeta_z\bar\zeta_w|
-{1\over\epsilon_\c'}\cdot{2\over\pi}\arctan\sqrt{\delta-1\over2}
-\epsilon_\c' r_\delta\biggr).}
$$
We now note that
$$
\eqalign{
\gr_{\bar D_\c(\epsilon_\c')}(z,w)&=\#(\Gamma\cap\{\pm1\})\cdot{1\over2\pi}
\log\biggl|{\zeta_z-\zeta_w\over 1-\zeta_z\bar\zeta_w}\biggr|\cr
&=\#(\Gamma\cap\{\pm1\})\left[
{1\over2\pi}\log|q_\c(z)-q_\c(w)|+{1\over\epsilon_\c'}
-{1\over2\pi}\log|1-\zeta_z\bar\zeta_w|\right].}
$$
Combining this with~\eqref{eq:grGamma-K} and the above bounds
on~$K(z,w)$ yields the proposition.\endproof

\appendixtrue

\unnumberedsection Appendix: Bounds on Legendre functions

In this appendix, we prove a number of bounds on the Legendre
functions $P^\mu_\nu(z)$ and~$Q^\mu_\nu(z)$ that are used in the rest
of the paper.

\proclaim Lemma. Let $u\in[1,3]$ and $\lambda\ge0$ be such that
$\lambda(u-1)\le{1\over2}$, and let $s\in\C$ be such that
$s(1-s)=\lambda$.  Then the real number $P^{-1}_{s-1}(u)$ satisfies
the inequalities
$$
(2-4/\pi)\sqrt{u-1\over u+1}
\le P^{-1}_{s-1}(u)
\le(4/\pi)\sqrt{u-1\over u+1}.
$$

\label{lemma:P-1-bound}

\proof We start by expressing the Legendre function~$P^\mu_\nu$ in
terms of Gau\ss's hypergeometric function~$F(a,b;c;z)$.  Because of
the many symmetries satisfied by the hypergeometric function (see
Erd\'elyi et~al.~\citex{Erdelyi}{Chapter~II}), there are lots of ways
to do this.  Using \citex{Erdelyi}{\S3.2, equation~3} gives
$$
P^{-1}_{s-1}(u)=\sqrt{u-1\over u+1}F\biggl(s,1-s;2;{1-u\over2}\biggr).
$$
Next we use the hypergeometric series for $F(a,b;c;z)$ with $z<1$
(see \citex{Erdelyi}{\S2.1, equation~2}):
$$
F(a,b;c;z)=\sum_{n=0}^\infty{(a)_n(b)_n\over(c)_n n!}z^n,
\eqnumber{hypergeometric-series}
$$
where
$$
(y)_n=\Gamma(y+n)/\Gamma(y)=y(y+1)\cdots(y+n-1).
$$
Putting $x={u-1\over2}$ for a moment,
using \eqref{hypergeometric-series} and applying the triangle
inequality, we get the bound
$$
\bigl|F(s,1-s;2;-x)-1\bigr|\le
\sum_{n\ge1}\left|{(s)_n(1-s)_n\over(2)_n n!}(-x)^n\right|
$$
The assumption $\lambda(u-1)\le{1\over2}$ is equivalent to $\lambda
x\le{1\over4}$.  Therefore the $n$-th term in the series on the
right-hand side can be bounded as follows:
$$
\eqalign{
\left|{(s)_n(1-s)_n\over(2)_n n!}(-x)^n\right|&=
{\prod_{k=0}^{n-1}\bigl(s(1-s)x+k(k+1)x\bigr)\over(2)_n n!}\cr
&\le{\prod_{k=0}^{n-1}\bigl({1\over4}+k(k+1)\bigr)\over(2)_n n!}\cr
&={\bigl({1\over2}\bigr)_n\bigl({1\over2}\bigr)_n\over(2)_n n!}.}
$$
This implies that
$$
\eqalign{
\bigl|F(s,1-s;2;-x\bigr)-1\bigr|
&\le F(\textfrac1/2,\textfrac1/2;2;1\bigr)-1\cr
&=4/\pi-1,}
$$
where the last equality follows from the formula
$$
F(a,b;c;1)={\Gamma(c)\Gamma(c-b-a)\over\Gamma(c-a)\Gamma(c-b)}
\quad\hbox{for $\Re c>0$ and $\Re c>\Re(a+b)$}
\eqnumber{eq:Fabc1}
$$
(see Erd\'elyi et~al.~\citex{Erdelyi}{\S2.1.3, equation~14} or
Iwaniec~\citex{Iwaniec}{equation~B.20}) and the fact that
$\Gamma(3/2)=\sqrt{\pi}/2$.  We conclude that
$$
\eqalign{
\biggl|P^{-1}_{s-1}-\sqrt{u-1\over u+1}\biggr|
&=\sqrt{u-1\over u+1}\bigl|F(s,1-s;2;-x\bigr)-1\bigr|\cr
&\le \sqrt{u-1\over u+1}(4/\pi-1)\cr
&=(4/\pi-1)\sqrt{u-1\over u+1},}
$$
which is equivalent to the inequalities in the statement of the
lemma.\endproof

\proclaim Lemma. For all real numbers $\nu\ge0$ and~$u>1$, the real
number $Q_\nu'(u)$ satisfies
$$
\biggl({2\over u+1}\biggr)^\nu{1\over u^2-1}
\le Q_\nu'(u)\le 0.
$$

\label{lemma:Q'-bound}

\proof We express $Q_\nu'$ in terms of the hypergeometric function
using \citex{Erdelyi}{\S3.6.1, equation~5, and \S3.2, equation~36}:
$$
Q_\nu'(u)=-\biggl({2\over u+1}\biggr)^\nu{1\over u^2-1}
{\Gamma(1+\nu)\Gamma(2+\nu)\over\Gamma(2+2\nu)}
F\biggl(\nu,1+\nu;2+2\nu;{2\over u+1}\biggr).
$$
Since $2/(u+1)<1$, the hypergeometric function is given by the
series~\eqref{hypergeometric-series}.  The non-negativity of all the
arguments gives the bounds
$$
\eqalign{
0&\le F\biggl(\nu,1+\nu;2+2\nu;{2\over u+1}\biggr)\cr
&\le F(\nu,1+\nu;2+2\nu;1)\cr
&={\Gamma(2+2\nu)\Gamma(1)\over\Gamma(2+\nu)\Gamma(1+\nu)};}
$$
the last equality follows from~\eqref{eq:Fabc1}.  Combining this with
the above formula for~$Q_\nu'(u)$ yields the claim.\endproof

\proclaim Lemma. For all real numbers $a$, $b$ and~$y$ with $b\ge
a>0$, we have
$$
\eqalign{
\exp\left(-{1\over 12}\left({1\over a}-{1\over b}\right)\right)
{(a^2+y^2)^{a/2-1/4}\over(b^2+y^2)^{b/2-1/4}}
&\le\left|\Gamma(a+iy)\over\Gamma(b+iy)\right|\cr
&\le\exp\left(b-a+{1\over 12}\left({1\over a}-{1\over b}\right)\right)
{(a^2+y^2)^{a/2-1/4}\over(b^2+y^2)^{b/2-1/4}}.}
$$

\label{lemma:gammaquot}

\proof We use Binet's formula for $\log\Gamma$ (see Erd\'elyi
et~al.\ \citex{Erdelyi}{\S1.9, equation~4}):
$$
\log\Gamma(z)=\bigl(z-\textfrac1/2\bigr)\log z-z+\textfrac1/2\log2\pi
+\int_0^\infty B(t)\exp(-zt)dt
\quad\hbox{for }\Re z>0,
$$
where
$$
\eqalign{
B(t)&=
\displaystyle{(\exp(t)-1)^{-1}-t^{-1}+\textfrac1/2\over t}\cr
&={{t\over2}\coth{t\over2}-1\over t^2}.}
$$
We write
$$
\eqalign{
M(a,y)&=\Re\int_0^\infty B(t)\exp(-(a+iy)t)dt\cr
&=\int_0^\infty B(t)\exp(-at)\cos(yt)dt.}
$$
Then we have
$$
\eqalign{
\log\left|{\Gamma(a+iy)\over\Gamma(b+iy)}\right|
&=\Re\Bigl((a+iy-1/2)\log(a+iy)-(b+iy-1/2)\log(b+iy)\Bigr)\cr
&\qquad-a+b+M(a,y)-M(b,y)\cr
&=(a-1/2)\log|a+iy|-(b-1/2)\log|b+iy|\cr
&\qquad-y\arg(a+iy)+y\arg(b+iy)-a+b+M(a,y)-M(b,y).}
$$
We note that
$$
y\arg(a+iy)-y\arg(b+iy)=y\arctan{(b-a)y\over ab+y^2}
\in[0,b-a].
$$
Using $b\ge a$, we conclude
$$
\log\left|{\Gamma(a+iy)\over\Gamma(b+iy)}\right|
\le(a-1/2)\log|a+iy|-(b-1/2)\log|b+iy|-a+b+M(a,y)-M(b,y)
$$
and
$$
\log\left|{\Gamma(a+iy)\over\Gamma(b+iy)}\right|
\ge(a-1/2)\log|a+iy|-(b-1/2)\log|b+iy|+M(a,y)-M(b,y).
$$

It remains to bound $M(a,y)-M(b,y)$.  The function $B$ satisfies
$$
0<B(t)\le\lim_{x\to 0} B(x)=1/12\quad\hbox{for all }t>0.
$$
Using this and the positivity of $\exp(-at)-\exp(-bt)$, we bound
$M(a,y)-M(b,y)$ as follows:
$$
\eqalign{
\bigl|M(a,y)-M(b,y)\bigr|
&\le\int_0^\infty B(t)\bigl(\exp(-at)-\exp(-bt)\bigr)
\left|\cos(yt)\right|dt\cr
&=\int_0^\infty B(t)\bigl(\exp(-at)-\exp(-bt)\bigr)dt\cr
&\le{1\over12}\int_0^\infty\bigl(\exp(-at)-\exp(-bt)\bigr)dt\cr
&={1\over12}\left({1\over a}-{1\over b}\right).}
$$
This implies the inequality we wanted to prove.\endproof


\proclaim Corollary.  For $b\ge a>0$, $b\ge1/2$, and $y\in\R$, we have
$$
\left|\Gamma(a+iy)\over\Gamma(b+iy)\right|\le
\exp\left(b-a+{1\over 12}\left({1\over a}-{1\over b}\right)\right)(a^2+y^2)^{-(b-a)/2}.
$$

\label{cor:gammaquot}





Given a real number $\sigma\in(0,1/2)$, we consider the strip
$$
S_\sigma=\{s\in\C\mid\sigma\le\Re s\le 1-\sigma\}.
$$
We put
$$
\eqalign{
C_\sigma&=\max\{1,\tan\pi\sigma\}(\sigma^{-1}-1)^{1/4}
\exp\left({{1\over2}+{1\over24\sigma(\textfrac1/2+\sigma)}}\right),\cr
C_\sigma'&=\max\{1,\tan\pi\sigma\}(\sigma^{-1}-1)^{1/4}
\exp\left({{1\over2}+{1\over24(1-\sigma)(\textfrac3/2-\sigma)}}\right).}
\eqnumber{eq:Csigma}
$$

\proclaim Proposition.  Let $m$ be an even non-negative integer, and
let $\sigma\in(0,1/2)$.  For all $s\in S_\sigma$ and all $u>1$, we
have
$$
\bigl|P^{-m}_{s-1}(u)\bigr|\le\bigl|s(1-s)\bigr|^{-(2m+1)/4}
{C_\sigma x^{m-\sigma}+C_\sigma' x^{m-1+\sigma}\over\sqrt{\pi}(4(u^2-1))^{m/2}}
\sum_{n=0}^\infty{|(\textfrac1/2-m)_n|\over n!}x^{-2n},
$$
where
$$
x=u+\sqrt{u^2-1},\quad u={x+x^{-1}\over2}.
$$

\label{prop:bound-Pm-regular}

\proof We use the following expression for $P^{-m}_{s-1}$ (see
Erd\'elyi et al.\ \citex{Erdelyi}{\S3.2, equation 27}):
$$
\eqalign{
P^{-m}_{s-1}(u)&=
{\Gamma(-\textfrac1/2+s)\over\sqrt\pi\Gamma(m+s\bigr)}
{x^{m+s-1}\over(x-x^{-1})^m}
F\bigl(\textfrac1/2-m,1-m-s;\textfrac3/2-s;x^{-2}\bigr)\cr
&\qquad+{\Gamma(\textfrac1/2-s)\over\sqrt\pi\Gamma(m+1-s\bigr)}
{x^{m-s}\over(x-x^{-1})^m}
F\bigl(\textfrac1/2-m,-m+s;\textfrac1/2+s;x^{-2}\bigr).}
$$
Using the hypergeometric series \eqref{hypergeometric-series} and the
functional equation
$$
\Gamma(z)\Gamma(1-z)={\pi\over\sin\pi z},
$$
we get
$$
\eqalign{
\sqrt\pi(x-x^{-1})^m P^{-m}_{s-1}(u)
&={\Gamma(-\textfrac1/2+s)\over\Gamma(m+s)}x^{m-1+s}
\sum_{n=0}^\infty{(\textfrac1/2-m)_n(-m+1-s)_n\over(\textfrac3/2-s)_n}
{x^{-2n}\over n!}\cr
&\qquad+{\Gamma(\textfrac1/2-s)\over\Gamma(m+1-s)}x^{m-s}
\sum_{n=0}^\infty{(\textfrac1/2-m)_n(-m+s)_n\over(\textfrac1/2+s)_n}
{x^{-2n}\over n!}\cr
&=\sum_{n=0}^\infty{(\textfrac1/2-m)_n\over n!}
\biggl\{{\Gamma(-\textfrac1/2+s)\over\Gamma(m+s)}
{\Gamma(\textfrac3/2-s)\Gamma(n-m+1-s)
\over\Gamma(n+\textfrac3/2-s)\Gamma(-m+1-s)}x^{m-1+s-2n}\cr
&\hskip2.8cm+{\Gamma(\textfrac1/2-s)\over\Gamma(m+1-s)}
{\Gamma(\textfrac1/2+s)\Gamma(n-m+s)
\over\Gamma(n+\textfrac1/2+s)\Gamma(-m+s)}x^{m-s-2n}\biggr\}\cr
&=\sum_{n=0}^\infty{(\textfrac1/2-m)_n\over n!}
\biggl\{{\sin\pi(m+s)\over\sin\pi(-\textfrac1/2+s)}
{\Gamma(n-m+1-s)\over\Gamma(n+\textfrac3/2-s)}x^{m-1+s-2n}\cr
&\hskip2.8cm+{\sin\pi(m+1-s)\over\sin\pi(\textfrac1/2-s)}
{\Gamma(n-m+s)\over\Gamma(n+\textfrac1/2+s)}x^{m-s-2n}\biggr\}.}
$$
Basic trigonometric manipulations simplify this to
$$
\eqalign{
P^{-m}_{s-1}(u)&={(-1)^m\tan(\pi s)\over\sqrt\pi(x-x^{-1})^m}
\sum_{n=0}^\infty{(\textfrac1/2-m)_n\over n!}
\biggl\{{\Gamma(n-m+s)\over\Gamma(n+\textfrac1/2+s)}x^{m-s-2n}\cr
&\hskip5.2cm-{\Gamma(n-m+1-s)\over\Gamma(n+\textfrac3/2-s)}
x^{m-1+s-2n}\biggr\}.}
\eqnumber{eq:nice-formula-P}
$$
On the right-hand side, the pole of $\tan(\pi s)$ at $s=1/2$ is
cancelled by a corresponding zero of the function defined by the sum.

For fixed $u>1$, we consider the holomorphic function
$$
\eqalign{
H(s)&=\bigl(s(1-s)\bigr)^{(2m+1)/4}\tan(\pi s)
\sum_{n=0}^\infty{(\textfrac1/2-m)_n\over n!}
\biggl\{{\Gamma(n-m+s)\over\Gamma(n+\textfrac1/2+s)}x^{m-s-2n}
\cr
&\hskip7cm-{\Gamma(n-m+1-s)\over\Gamma(n+\textfrac3/2-s)}
x^{m-1+s-2n}\biggr\}}
\eqnumber{eq:def-H}
$$
on~$S_\sigma$, where we have fixed a branch of
$s\mapsto(s(1-s))^{(2m+1)/4}$.  Because $H(s)=H(1-s)$, the
Phragm\'en--Lindel\"of principle gives
$$
\sup_{s\in S_\sigma}|H(s)|\le\sup_{y\in\R}|H(\sigma+iy)|.
$$
Together with \eqref{eq:nice-formula-P}, this implies
$$
|P^{-m}_{s-1}(u)|\le{\left|s(1-s)\right|^{-(2m+1)/4}
\sup_{y\in\R}|H(\sigma+iy)|\over\sqrt\pi(4(u^2-1))^{m/2}}
\quad\hbox{for all }s\in S_\sigma.
$$

Let $y\in\R$ and $s=\sigma+iy$.  Then we have
$$
\bigl|s(1-s)\bigr|^{(2m+1)/4}
=(\sigma^2+y^2)^{(2m+1)/8}((1-\sigma)^2+y^2)^{(2m+1)/8}.
$$
A straightforward calculation gives
$$
\left|\tan\pi s\right|=\left|\tan\pi(\sigma+iy)\right|\le
\max\{1,\tan\pi\sigma\}.
$$
Using Corollary~\ref{lemma:gammaquot} and the assumption that $m$ is
even, we bound the quotients of $\Gamma$-functions appearing on the
right-hand side of~\eqref{eq:def-H} independently of~$n$:
$$
\eqalign{
\left|{\Gamma(n-m+s)\over\Gamma(n+\textfrac1/2+s)}\right|
&={|\Gamma(n+s)|\over|\Gamma(n+\textfrac1/2+s)|}
\prod_{j=n-m}^{n-1}{1\over|j+s|}\cr
&\le\exp\left({1\over2}+{1\over24(n+\sigma)(n+\textfrac1/2+\sigma)}\right)
|n+\sigma+iy|^{-1/2}
{1\over|\sigma+iy|^{m/2}
|1-\sigma+iy|^{m/2}}\cr
&\le\exp\left({1\over2}+{1\over24\sigma(\textfrac1/2+\sigma)}\right)
{1\over(\sigma^2+y^2)^{(m+1)/4}((1-\sigma)^2+y^2)^{m/4}}.}
$$
This implies
$$
\eqalign{
\bigl|s(1-s)\bigr|^{(2m+1)/4}
\left|{\Gamma(n-m+s)\over\Gamma(n+\textfrac1/2+s)}\right|
&\le\exp\left({1\over2}+{1\over24\sigma(\textfrac1/2+\sigma)}\right)\cr
&\qquad\cdot{(\sigma^2+y^2)^{(2m+1)/8}((1-\sigma)^2+y^2)^{(2m+1)/8}\over
(\sigma^2+y^2)^{(m+1)/4}((1-\sigma)^2+y^2)^{m/4}}\cr
&=\exp\left({1\over2}+{1\over24\sigma(\textfrac1/2+\sigma)}\right)
{((1-\sigma)^2+y^2)^{1/8}\over(\sigma^2+y^2)^{1/8}}\cr
&\le\exp\left({1\over2}+{1\over24\sigma(\textfrac1/2+\sigma)}\right)
{(1-\sigma)^{1/4}\over\sigma^{1/4}}}
$$
and hence
$$
\bigl|s(1-s)\bigr|^{(2m+1)/4}|{\tan\pi s}|
\left|{\Gamma(n-m+s)\over\Gamma(n+\textfrac1/2+s)}\right|
\le C_\sigma.
$$
Similarly,
$$
\left|{\Gamma(n-m+1-s)\over\Gamma(n+\textfrac3/2-s)}\right|
\le\exp\left({1\over2}+{1\over24(1-\sigma)(\textfrac3/2-\sigma)}\right)
{1\over((1-\sigma)^2+y^2)^{(m+1)/4}(\sigma^2+y^2)^{m/4}}
$$
and
$$
\eqalign{
\bigl|s(1-s)\bigr|^{(2m+1)/4}
\left|{\Gamma(n-m+1-s)\over\Gamma(n+\textfrac3/2-s)}\right|
&\le\exp\left({1\over2}+{1\over24(1-\sigma)(\textfrac3/2-\sigma)}\right)\cr
&\qquad\cdot{(\sigma^2+y^2)^{(2m+1)/8}((1-\sigma)^2+y^2)^{(2m+1)/8}\over
((1-\sigma)^2+y^2)^{m/4}(\sigma^2+y^2)^{(m+1)/4}}\cr
&=\exp\left({1\over2}+{1\over24(1-\sigma)(\textfrac3/2-\sigma)}\right)
{((1-\sigma)^2+y^2)^{1/8}\over(\sigma^2+y^2)^{1/8}}\cr
&\le\exp\left({1\over2}+{1\over24(1-\sigma)(\textfrac3/2-\sigma)}\right)
{(1-\sigma)^{1/4}\over\sigma^{1/4}}.}
$$
This implies
$$
\bigl|s(1-s)\bigr|^{(2m+1)/4}|{\tan\pi s}|
\left|{\Gamma(n-m+1-s)\over\Gamma(n+\textfrac3/2-s)}\right|
\le C_\sigma'.
$$
We conclude that
$$
\eqalign{
\sup_{y\in\R}|H(\sigma+iy)|&\le
\sum_{n=0}^\infty{|(\textfrac1/2-m)_n|\over n!}
\bigl(C_\sigma x^{m-\sigma-2n}+C_\sigma' x^{m-1+\sigma-2n}\bigr)\cr
&=\bigl(C_\sigma x^{m-\sigma}+C_\sigma' x^{m-1+\sigma}\bigr)
\sum_{n=0}^\infty{|(\textfrac1/2-m)_n|\over n!}x^{-2n}.}
$$
This finishes the proof.\endproof





With $C_\sigma$ and~$C_\sigma'$ as in~\eqref{eq:Csigma}, we define an
elementary function $p_\sigma\colon[1,\infty)\to\R$ by
$$
p_\sigma(u)={C_\sigma x^{2-\sigma} + C_\sigma' x^{1+\sigma}\over 4\sqrt\pi}
\bigl((1-x^{-2})^{3/2}+3x^{-2}\bigr),
\quad\hbox{where }x=u+\sqrt{u^2-1}.
\eqnumber{eq:psig}
$$

\proclaim Corollary. For all $\sigma\in(0,1/2)$, $s\in S_\sigma$ and
$u>1$, we have
$$
|P^{-2}_{s-1}(u)|\le |s(1-s)|^{-5/4}{p_\sigma(u)\over u^2-1}.
$$

\label{cor:bound-P2-regular}

\proof We note that
$$
|(-\textfrac3/2)_n|=\cases{
-(-\textfrac3/2)_n& for $n=1$,\cr
(-\textfrac3/2)_n& otherwise.}
$$
This implies that for $z\in(0,1)$, we have
$$
\eqalign{
\sum_{n=0}^\infty{|(-\textfrac3/2)_n|\over n!}z^n
&=\sum_{n=0}^\infty{(-\textfrac3/2)_n\over n!}z^n
-2{(-\textfrac3/2)_1\over 1!}z^1\cr
&=(1-z)^{3/2}+3z.}
$$
The claim immediately follows from this identity and
Proposition~\ref{prop:bound-Pm-regular}.\endproof

\unnumberedsection References

\hyphenation{ma-the-ma-ti-sche}


\def\PSL{\mathop{\rm PSL}\nolimits}

\normalparindent=25pt
\parindent=\normalparindent
\parskip=1ex plus 0.5ex minus 0.2ex

\reference{Arakelov} {\cyr\S}.~{\cyr\Yu}.~{\cycsc \A\r\a\k\ye\l\o\v},
{\cyr \T\ye\o\r\i\ya\ \p\ye\r\ye\s\ye\ch\ye\n\i\j\ \d\i\v\i\z\o\r\o\v\ \n\a\
\a\r\i\f\m\ye\t\i\ch\ye\s\k\o\j\ \p\o\v\ye\r\kh\n\o\s\t\i}.\hfil\break
{\cyti \I\z\v\ye\s\t\i\ya\ \A\k\a\d\ye\m\i\i\ \N\a\u\k\ \S\S\S\R{\rm
  ,} \s\ye\r\i\ya\ \m\a\t\ye\m\a\t\i\ch\ye\s\k\a\ya\/} {\bf 38}
  (1974), {\cmcyr\char"19}~6, 1179--1192.
\hfil\break
S. Yu.\ {\sc Arakelov}, Intersection theory of divisors on an
arithmetic surface.  {\it Mathematics of the USSR Izvestiya\/} {\bf 8}
(1974), 1167--1180.  (English translation.)


\reference{thesis} P. J. {\sc Bruin}, {\sl Modular curves, Arakelov
theory, algorithmic applications\/}.  Proefschrift (Ph.\thinspace D.\
thesis), Universiteit Leiden, 2010.

\reference{book} S. J. {\sc Edixhoven} and J.-M. {\sc Couveignes}
(with R. S. {\sc de Jong}, F. {\sc Merkl} and J. G. {\sc Bosman}),
{\sl Computational Aspects of Modular Forms and Galois
Representations\/}.  Annals of Mathematics Studies {\bf 176}.
Princeton University Press, 2011.

\reference{Erdelyi} A. {\sc Erd\'elyi}, W. {\sc Magnus}, F. {\sc
Oberhettinger} and F. G. {\sc Tricomi}, {\sl Higher Transcendental
Functions\/}, Volume~I\null.  Bateman Manuscript Project, California
Institute of Technology.  McGraw-Hill, New York/Toronto/London, 1953.


\reference{Faddeev} {\cyr\L}.~{\cyr\D}.~{\cycsc \F\a\d\d\ye\ye\v},
{\cyr \R\a\z\l\o\zh\ye\n\i\ye\ \p\o\ \s\o\b\s\t\v\ye\n\n\y\m\
\f\u\n\k\ts\i\ya\m\ \o\p\ye\r\a\t\o\r\a\ \L\a\p\l\a\s\a\ \n\a\
\f\u\n\d\a}-{\cyr\m\ye\n\t\a\l\soft\n\o\j\ \o\b\l\a\s\t\i\ \d\i\s\k\r\ye\t\o\j\
\g\r\u\p\p\y\ \n\a\ \p\l\o\s\k\o\s\t\i\ \L\o\b\a\ch\ye\v\s\k\o\g\o}.
{\cyti \T\r\u\d\y\ \M\o\s\k\o\v\s\k\o\g\o\
\M\a\t\ye\m\a\t\i\ch\ye\s\k\o\g\o\ \O\b\shch\ye\s\t\v\a}
{\bf 17} (1967), 323--350.\hfill\break
L. D. {\sc Faddeev}, Expansion in eigenfunctions
of the Laplace operator on the fundamental domain of a discrete group
on the Loba\v{c}evski\u\i\ plane.  {\it Transactions of the Moscow
Mathematical Society\/} {\bf 17} (1967), 357--386.  (English
translation.)

\reference{Faltings: Calculus on arithmetic surfaces} G. {\sc
Faltings}, Calculus on arithmetic surfaces.  {\it Annals of
Mathematics\/}~(2) {\bf 119} (1984), 387--424.

\reference{Fay} J. D. {\sc Fay}, Fourier coefficients of the resolvent
for a Fuchsian group.  {\it Journal f\"ur die reine und angewandte
Mathematik\/} {\bf 293/294} (1977), 143--203.



\reference{Hejhal-II} D. A. {\sc Hejhal}, {\sl The Selberg trace
formula for~$\PSL(2,{\bf R})$\/}, Volume~2.  Lecture Notes in
Mathematics~{\bf 1001}.  Springer-Verlag, Berlin/Heidelberg, 1983.

\reference{Huber} H. {\sc Huber}, \"Uber eine neue Klasse automorpher
Funktionen und ein Gitterpunktproblem in der hyperbolischen
Ebene.\ I\null.  {\it Commentarii Mathematici Helvetici\/} {\bf 30}
(1956), 20--62.

\reference{Iwaniec} H. {\sc Iwaniec}, {\sl Introduction to the
Spectral Theory of Automorphic Forms\/}.  Revista Mate\-m\'a\-tica
Ibero\-americana, Madrid, 1995.

\reference{Jorgenson-Kramer: Green's functions} J. {\sc Jorgenson} and
J. {\sc Kramer}, Bounds on canonical Green's functions.  {\it
Compositio Mathematica\/} {\bf 142} (2006), no.~3, 679--700.

\reference{Kim} H. H. {\sc Kim}, Functoriality for the exterior square
of~${\rm GL}_4$ and the symmetric fourth of~${\rm GL}_2$.  With
appendix~1 by D. {\sc Ramakrishnan} and appendix~2 by {\sc Kim} and
P.~{\sc Sarnak}.  {\it Journal of the A.M.S.\/} {\bf 16} (2002),
no.~1, 139--183.


\reference{Maass} H. {\sc Maa\ss}, \"Uber eine neue Art von
nichtanalytischen automorphen Funktionen und die Be\-stim\-mung
Dirichletscher Reihen durch Funktionalgleichungen.  {\it Mathematische
Annalen\/} {\bf 121} (1949), 141--183.



\reference{pari} PARI/GP, version~2.5.1.  Bordeaux, 2011, {\tt
  http://pari.math.u-bordeaux.fr/}.

\reference{Patterson} S. J. {\sc Patterson}, A lattice problem in
hyperbolic space.  {\it Mathematika\/} {\bf 22} (1975), 81--88.

\reference{Selberg: Harmonic analysis} A. {\sc Selberg}, Harmonic
analysis and discontinuous groups in weakly symmetric Riemannian
spaces with applications to Dirichlet series.  {\it Journal of the
Indian Mathematical Society (N.S.)\/}\ {\bf 20} (1956), 47--87.
(=~{\sl Collected Papers\/}, Volume~1, 423--463.  Springer-Verlag,
Berlin, 1989.)

\reference{Selberg: ICM} A. {\sc Selberg}, Discontinuous groups and
harmonic analysis.  In: {\sl Proceedings of the International Congress
of Mathematicians (Stockholm, 15--22 August 1962)\/}, 177-189.
Institut Mittag-Leffler, Djursholm, 1963.  (=~{\sl Collected Papers\/},
Volume~1, 493--505.  Springer-Verlag, Berlin, 1989.)

\reference{Selberg: Fourier coefficients} A. {\sc Selberg}, On the
estimation of Fourier coefficients of modular forms.  In: A.~L. {\sc
Whiteman} (editor), {\sl Theory of Numbers\/}, 1--15.  Proceedings of
Symposia in Pure Mathematics, VIII\null.  American Mathematical
Society, Providence, RI, 1965.  (=~{\sl Collected Papers\/}, Volume~1,
506--520.  Springer-Verlag, Berlin, 1989.)



\vskip2.5cm

\leftline{Peter Bruin}
\leftline{Institut f\"ur Mathematik}
\leftline{Universit\"at Z\"urich}
\leftline{Winterthurerstrasse 190}
\leftline{CH-8057 Z\"urich}
\smallskip
\leftline{\tt peter.bruin@math.uzh.ch}

\bye